\documentstyle[amscd,amssymb,verbatim,11pt]{amsart}

\theoremstyle{plain}
\newtheorem{Prop}{Proposition}[section]
\newtheorem{Thm}[Prop]{Theorem}
\newtheorem{Cor}[Prop]{Corollary}

\newtheorem{Lem}[Prop]{Lemma}
\newtheorem{notation}[Prop]{Notation}

\theoremstyle{definition}
\newtheorem{Def}[Prop]{Definition}

\theoremstyle{remark}
\newtheorem{Rem}[Prop]{Remark}

\def\dim{\mathop{\roman{dim}}}
\def\int{\mathop{\roman{int}}}

\def\1{^{-1}}
\def\Tilde{\widetilde}
\def\Z{{\bold Z}}
\def\Q{{\bold Q}}

\def\p{{\bold p}}
\def\q{{\bold q}}
\def\B{\cal B_{\cal G}}

\def\im{\text{Im}}
\def\dim{\text{dim}}

\def\Tor{\text{Tor}}

\def\ExD{\text{ext--dim}}

\def\int{\text{Int}}

\numberwithin{equation}{section}
\begin{document}

\title[
Algebra of dimension theory
]%
   {
Algebra of dimension theory
}

\author{
Jerzy Dydak
}
\date{ April 4, 2004
}


\thanks{ Research supported in part by a grant
 DMS-0072356 from the National Science Foundation
}

\begin{abstract}
The dimension algebra of graded groups is introduced. With the help of known geometric results of
extension theory that algebra induces all known results of the cohomological
dimension theory. Elements of the algebra are equivalence
classes $\dim(A)$ of graded groups $A$. There are two geometric interpretations
of those equivalence classes:
\linebreak 1. For pointed CW complexes $K$ and $L$,
$\dim(H_\ast(K))=\dim(H_\ast(L))$ if and only if
the infinite symmetric products $SP(K)$ and $SP(L)$ are of the same
extension type (i.e., $SP(K)\in AE(X)$ iff $SP(L)\in AE(X)$ for all compact $X$).
\linebreak 2. For pointed compact spaces $X$ and $Y$,
$\dim(\cal H^{-\ast}(X))=\dim(\cal H^{-\ast}(Y))$
if and only if $X$ and $Y$ are of the same dimension type
(i.e., $\dim_G(X)=\dim_G(Y)$ for all Abelian groups $G$).
\par Dranishnikov's version of Hurewicz Theorem in extension theory
becomes
$\dim(\pi_\ast(K))=\dim(H_\ast(K))$ for all simply connected $K$.
\par The concept of cohomological dimension $\dim_A(X)$
of a pointed compact space $X$ with respect to a graded group $A$ is introduced.
It turns out $\dim_A(X)\leq 0$ iff $\dim_{A(n)}(X)\leq n$
for all $n\in\Z$. If $A$ and $B$ are two positive graded groups,
then $\dim(A)=\dim(B)$ if and only if $\dim_A(X)=\dim_B(X)$
for all compact $X$.
\end{abstract}

\maketitle

\medskip
\medskip
\tableofcontents

\section{Introduction }

\begin{notation} \label{XX1.0} Throughout the paper $K$, $L$, and $M$ are
reserved for pointed CW complexes. 
$SP(K)$ (see \cite{Git}, p.168) is the
{\bf infinite symmetric product} of a pointed CW complex $K$.
$X$ and $Y$ are general topological spaces
(quite often compact or compact metrizable).
We will frequently omit
coefficients in the case of
integral homology and cohomology. Thus,
 $H_n(K;\Z)$ will be shortened to $H_n(K)$ and $H^n(X;\Z)$ will be
shortened to $H^n(X)$.
\end{notation} 

In the paper we consider {\bf graded groups} $A$ indexed by integers.
The $n$-th term of $A$ will be denoted by $A(n)$, so $A=\{A(n)\}_{n\in\Z}$.
\par Here are the main examples of graded groups from the point of view of this paper.

\begin{Def} \label{XXi.2}
The {\bf homology graded group} $H_\ast(K)$ of a pointed CW complex
$K$ defined by $(H_\ast(K))(n)=H_n(K)$ for $n\in\Z$.
\par The {\bf reversed cohomology graded group} $H^{-\ast}(X)$
of a pointed compact space $X$ is defined by $(H^{-\ast}(X))(n)=H^{-n}(X)$
for $n\in\Z$.
\par The {\bf reversed total cohomology graded group} $\cal H^{-\ast}(X)$
of a pointed compact space $X$ is defined by declaring $(\cal H^{-\ast}(X))(n)$
to be the direct sum of all $H^{-n}(X,A)$, where $A$ ranges
over all pointed closed subsets of $X$ (the total cohomology group
was introduced by Shchepin \cite{S} without reversing indices).
\end{Def}

The main concept of the paper is the {\bf homological dimension}
$\dim_G(A)$ of a graded group $A$ with respect to an Abelian group $G$.
In the case of $A=H_\ast(K)$ it is equal to the homological
dimension $\dim_G(K)$ introduced in \cite{Dy$_6$} as the supremum
of $n$ such that $H_k(K;G)=0$ for all $k<n$.
In the case of $A=\cal H^{-\ast}(X)$ it is equal to the negative
of the cohomological dimension $\dim_G(X)$ of $X$ (it is the infimum
of $n$ such that $K(G,n)$ is the absolute extensor of $X$).
\par Using the homological dimension of graded groups we give a uniform
description of the Bockstein Theory, the Bockstein Algebra (see \cite{K}), and of
the Dual Bockstein Algebra (see \cite{D-D$_1$}). 

We will use the following results of geometric nature.

\begin{Thm} \label{XXXintro.1} (see \cite{D$_3$}) 
Suppose $K$ is a pointed connected CW complex. If $X$ is finite-dimensional and
$K(\pi_i(K),i)\in AE(X)$ for each $i\ge 1$, then $K\in AE(X)$.
\end{Thm} 

Of major importance to us is the following result of Dranishnikov:
\begin{Thm} \label{XXX2.10a} (see \cite{D$_5$}) 
Suppose $K$ is a pointed CW complex. If $X$ is compact and $K\in AE(X)$, then
$SP(K)\in AE(X)$. 
\end{Thm}

Since $SP(K)$ is homotopy equivalent
to the weak product of Eilenberg-MacLane spaces
$K(\tilde H_i(K),i)$ (see \cite{Git}, Corollary 6.4.17 on p.223)
one has the following:
\begin{Thm} \label{XXX2.10b} (see \cite{D$_5$}) 
Suppose $K$ is a pointed CW complex. If $X$ is compact, then
$SP(K)\in AE(X)$
is equivalent to $K(\tilde H_i(K),i)\in AE(X)$ for all $i\ge 0$. 
\end{Thm}

The author would like to thank Akira Koyama for discussions on the subject
of the paper, and E.Shchepin for giving a series of excellent talks
(in Russian) during the workshop \lq Algebraic ideas in dimension
theory\rq \ held in Warsaw (Fall 1998). Those talks were the starting
point in author's understanding of paper \cite{S}
which eventually led to the ideas developed in this paper.

\section{Algebras with valuations}
The basic
algebraic structure underlying most of algebraic topology
is that of a family of homeomorphism classes
of pointed spaces belonging to some specialized class closed under
 operations of wedge $X\vee Y$
and the smash product $X\wedge Y$. 
Both operations are commutative, 
associative, have neutral elements, and one has the distributivity
$X\wedge (Y\vee Z)=(X\wedge Y)\vee (X\wedge Z)$.
Let us call it the {\bf Standard Algebra}.
Notice that the wedge plays the role analogous
to the direct sum of vector spaces in the sense
that $Map(X\vee Y,Z)=Map(X,Z)\times Map (Y,Z)$,
and the smash product plays the role of the
tensor product of vector spaces
in the sense that $Map(X\wedge Y,Z)=Map(X,Map (Y,Z))$.
The formal categorical analogy with linear algebra can be worded
as follows: the wedge functor is adoint to the cartesian product
functor, and the smash product functor is adjoint to the
$hom$ functor sending $X$ to $Map(X,Z)$.
\par 
Suppose we have an algebra $\cal A$ and a {\bf discrete valuation}
 $v(X)$ associated with elements of $\cal A$, so that 
the following conditions are satisfied:
\par{1.} $v(X)\in \Z\cup\{\pm\infty\}$.
\par{2.} $v(X\vee Y)=\min(v(X),v(Y))$.
\par{3.} $v(X\wedge Y)\ge v(X)+v(Y)$.
\par
Using valuation $v$ one can establish a new equivalence
relation between elements of $\cal A$ and obtain its quotient
$\cal A/\sim$. Here is how it goes: $X\sim Y$
means $v(X\wedge Z)=v(Y\wedge Z)$ for all $Z\in \cal A$.
It is easy to check that $[X]\vee [Y]:=[X\vee Y]$
and $[X]\wedge [Y]:=[X\wedge Y]$ are well-defined.
Obviously, the valuation $v$ factors through $\cal A/\sim$.
Moreover, in the new algebra, $v(a\cdot c)=v(b\cdot c)$ for all $c$
implies $a=b$.
\par Let us point out similarity of our discrete valuations to those
used in classical algebra (see \cite{Krull1931}). 
There, $v$ is defined on $K\setminus\{0\}$, where
$K$ is a field and the following conditions are satisfied:
\par{1.} $v(a)\in \Z$.
\par{2.} $v(a+b)\ge \min(v(a),v(b))$.
\par{3.} $v(a\cdot b)=v(a)+v(b)$.
\par

Let us present three basic examples.
\par\noindent \centerline {\bf Kuzminov Algebra} 
$\cal A$ is the subalgebra  of the Standard Algebra
consisting of all  pointed metrizable compacta of finite dimension
and $v(X)$
is the negative of the covering dimension of $X$.
We will call the resulting quotient algebra $\cal A/\sim$
the {\bf Kuzminov Algebra}.

\par\noindent \centerline {\bf Shchepin Algebra}
\par $\cal A$ is the subalgebra  of the Standard Algebra
consisting of all pointed CW complexes and $v(K)$
is the {\bf connectivity index} $cin(K)$
of $K$, i.e. 
the supremum of $n$ so that $H_k(K)=0$ for all $k < n$.
We will call the resulting quotient algebra $\cal A/\sim$
the {\bf Shchepin Algebra} as its roots can be found in \cite{S}.

\par\noindent \centerline {\bf Dimension Algebra of Graded Groups}
\par The elements of $\cal A$
are isomorphism classes of graded groups $A$.
The wedge $A\vee B$ is simply the direct sum of $A$ and $B$, and
 the product $A\wedge B$ is defined so that:
\par 1. $X\to H^{-\ast}(X)$
is a homomorphism of the Standard Algebra of pointed compact spaces
to $\cal{A}$.
\par 2. $K\to H_\ast(K)$ is a homomorphism of the Standard Algebra of pointed CW complexes
to $\cal{A}$.

\par The valuation is the connectivity index $cin(A)$ of a graded
group $A$. It is the supremum of $n$ such
that $A_i=0$ for all $i<n$.
The quotient algebra is called the {\bf Dimension Algebra of Graded Groups}
(denoted by $\cal{DGG}$) and will be introduced in Section 6.

\section{Extension algebras}

In this section we introduce geometrically two important algebras
related to extension properties of pointed countable CW complexes.
\par Recall the following definition from \cite{D-D$_1$}.

\begin{Def}\label{XXea.1}
Let $\cal C$ be a class of compact spaces. Given two
pointed CW complexes $K$ and $L$, $K\sim_{\cal C} L$
means that the subclass of $\cal C$
consisting of all $X$ such that $K\in AE(X)$ is identical
the subclass of $\cal C$
consisting of all $X$ such that $L\in AE(X)$.
In particular, $K\sim_X L$ means $K\sim_{\{X\}} L$.
\end{Def}

The following result follows from Proposition 3.3 of \cite{D-D$_1$}.
However, since it is purely geometric we repeat its proof here.
\begin{Lem}\label{XXea.2}
Suppose $K_1,K_2,L_1,L_2$ are pointed countable CW complexes.
If $K_1\sim_X K_2$ and $L_1\sim_X L_2$ for all compacta
(respectively, all finite-dimensional compacta) $X$,
then $\Sigma(K_1\wedge L_1)\sim_X \Sigma(K_2\wedge L_2)$ for all compacta
(respectively, all finite-dimensional compacta) $X$.
\end{Lem}
\begin{pf} It suffices to show that $\Sigma(K_1\wedge L_1)\in AE(X)$
implies $\Sigma(K_2\wedge L_2)\in AE(X)$ for all compacta
(respectively, finite-dimensional compacta) $X$.
Since $\Sigma(K\wedge L)$ is homotopy equivalent to the join
$K\ast L$ of $K$ and $L$, $X$ can be expressed
as the union of two of its subsets $A$ and $B$ so that
$K_1\in AE(A)$, $L_1\in AE(B)$, and $A$ is a countable
union of closed subsets of $X$ (see \cite{D$_4$}).
Therefore $K_2\in AE(A)$ and, by Olszewski's Completion Theorem
(see \cite{O$_1$}) there is a $G_\delta$-subset $A'$ of $X$
containing $A$ such that $K_2\in AE(A')$.
Now, $X\setminus A'$ is an $F_\sigma$-subset of $X$ satisfying
$L_1\in AE(X\setminus A')$. Therefore $L_2\in AE(X\setminus A')$
and the main result of \cite{Dy$_3$} says that $K_2\ast L_2\in AE(X)$
which is the same as saying $\Sigma(K_2\wedge L_2)\in AE(X)$.
\end{pf}

\begin{Cor}\label{XXea.3}
Suppose $K_1,K_2,L_1,L_2$ are pointed countable CW complexes.
\par 1. If $K_1\sim_X L_1$ for all compacta
(respectively, all finite-dimensional compacta) $X$,
then $\Sigma(K_1)\sim_X \Sigma(L_1)$ for all compacta
(respectively, all finite-dimensional compacta) $X$.
\par 2. Suppose $\Sigma^m(K_1)\sim_X \Sigma^m(L_1)$ for some $m\ge 0$ and for all compacta
(respectively, all finite-dimensional compacta) $X$.
If  $\Sigma^n(K_2)\sim_X \Sigma^n(L_2)$ for some $n\ge 0$ and for all compacta
(respectively, all finite-dimensional compacta) $X$,
then $\Sigma^{m+n+1}(K_1\wedge L_1)\sim_X \Sigma^{m+n+1}(K_2\wedge L_2)$ for all compacta
(respectively, all finite-dimensional compacta) $X$.
\end{Cor}
\begin{pf} Put $L=S^0$. Since $K\wedge L=K$ for all $K$, 1) follows from
\ref{XXea.2}.
\par 2) also follows from \ref{XXea.2} as $\Sigma(\Sigma^m(K)\wedge \Sigma^n(L))$
is homotopy equivalent to $\Sigma^{m+n+1}(K\wedge L)$.
\end{pf}

Corollary \ref{XXea.3} means that the following definitions make sense.

\begin{Def}\label{XXea.4}
Consider the set of homeomorphic classes of all pointed countable CW complexes.
Two pointed CW complexes $K$ and $L$ are said to be equivalent
if there is $m\ge$ such that $\Sigma^m(K)\sim_X \Sigma^m(L)$ for all finite-dimensional compacta $X$.
The algebra with addition defined by $[K]+[L]=[K\vee L]$
and multiplication defined by $[K]\cdot [L]=[K\wedge L]$
will be called the {\bf Dranishnikov-Dydak Algebra}
as it is equivalent to the one introduced in  \cite{D-D$_1$}. 
\end{Def}

\begin{Def}\label{XXea.5}
Consider the set of homeomorphic classes of all pointed countable CW complexes.
Two pointed CW complexes $K$ and $L$ are said to be equivalent
if there is $m\ge$ such that $\Sigma^m(K)\sim_X \Sigma^m(L)$ for all compacta $X$.
The algebra with addition defined by $[K]+[L]=[K\vee L]$
and multiplication defined by $[K]\cdot [L]=[K\wedge L]$
will be called the {\bf Stable Extension Algebra}. 
\end{Def}

\begin{Rem}\label{XXea.6}
Notice that Dranishnikov-Dydak Algebra and the Stable Extension Algebra
are not identical. Indeed, $K(Z,n)\sim_X S^n$ for all $n\ge 1$
and all finite dimensional compacta $X$ (this amounts to the Alexandroff
Theorem stating that, in the class of finite-dimensional
compacta, both covering dimension and the integral cohomological dimension
coincide). However, for all $n\ge 2$ and all $m\ge 0$ there exist 
compacta $X$ such that $\Sigma^m(K(Z,n))\sim_X \Sigma^m(S^n)$
fails (see Theorem 9.3 of \cite{Dy$_6$}).
\end{Rem}

\section{Homological dimension of graded groups}

\begin{Def} \label{XXh.1} Two graded groups $A$ and $B$ are
{\bf isomorphic} (notation: $A\equiv B$) if $A(n)$ is isomorphic
to $B(n)$ for each $n\in\Z$.
\par Given a family $\{A_t\}_{t\in T}$ of graded groups,
its {\bf direct sum} $\bigoplus\limits_{t\in T} A_t$
is defined as the graded group whose $n$-th term
is $\bigoplus\limits_{t\in T} A_t(n)$.
\end{Def}

\par Our first step is to generalize the concept of smash product
($K\wedge L$ for pointed CW complexes and $X\wedge Y$ for pointed compact spaces)
to graded groups.


\begin{Def} \label{XX4.1} Given graded groups $A$ and $B$ their {\bf smash
product}
$A\wedge B$ is defined by
$(A\wedge B)(n)=\bigoplus \{A(k)\otimes B(n-k)\mid k\in\Z\}\oplus
\bigoplus \{A(k)*B(n-k-1)\mid k\in\Z\}$.

\end{Def} 

\begin{Def} \label{XX4.2} The {\bf suspension operator} $\Sigma^k$,
$k\in\Z\cup\{\pm\infty\}$,
on graded groups is defined as follows:
\par{1.} if $k$ is an integer, then $\Sigma^k(A)(n)=A(n-k)$ for all $n$.
\par{2.} if $k=-\infty$, then
$\Sigma^k(A)(n)=\bigoplus\limits_{m\in\Z}A(m)$ for all $n$.
\par{3.} if $k=\infty$, then $\Sigma^k(A)(n)=0$ for all $n$.

\end{Def}

\begin{Rem}  Notice that $\Sigma^1(H_*(K))=H_*(\Sigma(K))$ for pointed CW complexes
and $\Sigma^{-1}(H^{-\ast}(X))=H^{-\ast}(\Sigma(X))$ for pointed compact spaces.
 
\end{Rem}


\begin{notation} \label{XX4.3} When convenient, a group $G$ will be identified
with the graded group $A$ so that $A(0)=G$ and $A(k)=0$ if $k\ne 0$.
In particular, we will talk about graded groups $\Sigma^k(G)$
if $G$ is a group. One should think of groups as analogs of Moore spaces.

\end{notation}

\begin{Rem} In view of \ref{XX4.3}
the smash product $A\wedge G$ of $A$ and a group $G$
can be viewed as $(A\otimes G)\oplus \Sigma(A\ast G)$.
\end{Rem}

\begin{Def} \label{XX4.4} Suppose $A$ is a graded Abelian group and $G$
is an Abelian group.
The {\bf homological dimension} $\dim_G(A)\in\Z\cup\{\pm\infty\}$ of $A$
is defined as
$$\sup\{k\in\Z\cup\{\pm\infty\}\mid (A\wedge G)(n)=0\text{ for all }n<k\}.$$
\par $\dim(A)\leq\dim(B)$ means that $\dim_F(A)\leq\dim_F(B)$
for any Abelian group $F$.

\end{Def}


\begin{Prop} \label{XX4.5} 
\par{a.} For every graded group $A$ there is a free
chain complex $C$ such that $A=H_*(C)$.
\par{b.}
Suppose $A$ is the homology $H_*(C)$ of a free chain complex $C$,
$B$ is the homology $H_*(D)$ of a free chain complex $D$, and $G$ is an
Abelian group.
The following isomorphisms hold:
\par {1.} $A\wedge G\equiv H_*(C\otimes G)$,
\par {2.} $A\wedge B\equiv H_*(C\otimes D)$.

\end{Prop}

\begin{pf}  Given a graded group $A$ we can treat it as a trivial chain
complex
(i.e., all the boundary homomorphisms are trivial). Notice that
$A=H_*(A)$ in such case. Lemma 12 in \cite{Sp} (p.225) says that
there is a free approximation $C$ of $A$. In particular,
$H_*(C)=H_*(A)=A$.
\par
1) is a consequence of the Universal-Coefficient Theorem for homology
(see \cite{Sp}, Theorem 8 on p.222). 2) is a consequence of the K\" unneth Formula
for homology (see \cite{Sp}, Theorem 3 on p.230).
 \end{pf}


\begin{Prop} \label{XX4.7} Smash product is associative and commutative.

\end{Prop}

\begin{pf}  The commutativity of smash product follows from commutativity
of the tensor product and the torsion product.
Given three graded groups $A_i$, $1\leq i\leq 3$, choose
free chain complexes $C_i$ such that $A_i=H_*(C_i)$ for $1\leq i\leq 3$
(see a) of \ref{XX4.5}). Using \ref{XX4.5} notice that each of $(A_1\wedge A_2)\wedge A_3$
and $A_1\wedge (A_2\wedge A_3)$ equals $H_*(C_1\otimes C_2\otimes C_3)$
(use associativity of the tensor product).
 \end{pf}


\begin{Cor} \label{XX4.8} Suppose $A_1,A_2$ and $B_1,B_2$
are graded groups. If $\dim(A_i)\leq \dim(B_i)$ for $i\leq 2$,
then $\dim(A_1\oplus B_1)\leq \dim(A_2\oplus B_2)$
and
$\dim(A_1\wedge B_1)\leq \dim(A_2\wedge B_2)$.
\end{Cor} 
\begin{pf}  The inequality $\dim(A_1\oplus B_1)\leq \dim(A_2\oplus B_2)$
follows from $(A\oplus B)\wedge G\equiv (A\wedge G)\oplus (B\wedge G)$
for any graded groups $A$, $B$ and any group $G$.
\par First, consider the case $B_1=B_2=B$ is a group.
Notice that $A_1\leq A_2$
means $\dim_\Z(A_1\wedge G)\leq \dim_\Z(A_2\wedge G)$
for all groups $G$. Therefore
$\dim_\Z(A_1\wedge B)\leq \dim_\Z(A_2\wedge B)$
for all $B=\Sigma^n(G)$, where $n\in\Z$ and $G$ is a group.
Since every graded group $C$ is a direct sum of
$\Sigma^n(C(n))$, $n\in\Z$, one gets $\dim_\Z(A_1\wedge C)\leq \dim_\Z(A_2\wedge C)$
for all $C$. In particular, for $C=B\wedge G$
this amounts to $\dim_G(A_1\wedge B)\leq \dim_G(A_2\wedge B)$.
\par The general case follows from commutativity of the smash product
and the special case;
$\dim(A_1\wedge B_1)\leq \dim(A_1\wedge B_2)\leq \dim(A_2\wedge B_2)$.
 \end{pf}

For any two Abelian groups $F$ and $G$ one can discuss the concept
of $\dim(G)\leq \dim(F)$ and of $\dim(G)=\dim(F)$ using Convention \ref{XX4.3}.


\begin{Prop} \label{XX4.9} For any two Abelian groups $F$ and $G$, the
dimension
$\dim_G(F)$ can attain only three values; $0$, $1$, and $\infty$.
\par{1.} $\dim_G(F)=\infty$ if and only if $G\otimes F=0=G*F$.
\par{2.} $\dim_G(F)=1$ if and only if $G\otimes F=0$ and $G*F\ne 0$.
\par{3.} $\dim_G(F)=0$ if and only if $G\otimes F\ne 0$.

\end{Prop}

\begin{pf}  Let $A=F$ and $B=A\wedge G$. Notice that
$B(1)=F*G$, $B(0)=F\otimes G$, and $B(n)=0$ for $n\ne 0,1$.
 \end{pf}


\begin{Cor} \label{XX4.10} The following conditions are equivalent
for any two Abelian groups $F$ and $G$:
\par{1.} $\dim(F)\leq\dim(G)$.
\par{2.} $\dim_F(A)\leq\dim_G(A)$ for any graded group $A$.
\par{3.} If $H\otimes F=0$, then $H\otimes G=0$ for every Abelian group $H$.
If $H\otimes F=0$ and $H*F=0$, then $H\otimes G=0$ and $H*G=0$ for every
Abelian group $H$.

\end{Cor}

\begin{pf}  1)$\iff$ 2) is immediate from \ref{XX4.9}. 3)$\implies$ 2) is obvious
and 2)$\implies$ 3) by using $A=H$.
 \end{pf}

\begin{Rem}  Notice how these concepts relate to Shchepin's \cite{S}
factorial
domination.
 
\end{Rem}

\section{Extension theory and homological dimension}

\begin{Cor} \label{XXeh.1} Let $G$ be an Abelian group.
If $K$ and $L$ are pointed CW complexes, then
\par 1. $H_\ast(K;G)\equiv H_\ast(K)\wedge G$.
\par 2. $H_\ast(K\wedge L)\equiv H_\ast(K)\wedge H_\ast(L)$.
\end{Cor} 
\begin{pf}  For any CW complex $M$ one has the free chain complex
$C(M)$ so that $C(M)(n)=C_n(M)$ is the group of cellular $n$-chains.
$H_*(M)$ is simply $H_*(C(M))$ and $H_*(M;G)$ is defined as
$H_*(C(M)\otimes G)$. Thus, 1) of \ref{XXeh.1} follows from \ref{XX4.5}.
2) of \ref{XXeh.1} follows from the K\" unneth formula for singular homology
(see Theorem 10 on p.235 in \cite{Sp}).
 \end{pf}

In \cite{Dy$_6$} the inequality $K\leq_G L$ was defined to mean
$\dim_G(K)\leq \dim_G(L)$, where $\dim_G(M)$
is the minimum of $n$ such that $H_k(M;G)=0$ for all $k < n$.

\begin{Cor} \label{XX4.6} Suppose $K$ and $L$ are CW complexes and
$G$ is an Abelian group. $K\leq_GL$ if and only if $\dim_G(\tilde
H_*(K))\leq\dim_G(\tilde H_*(L))$.
\end{Cor}

Given a cohomology theory $h^\ast$, by $h^{-\ast}(X)$
we will denote {\bf the reversed cohomology graded group}
defined by $(h^{-\ast}(X))(n)=h^{-n}(X)$.

Obviously, every homology theory $h_*$
(respectively, cohomology theory $h^*$) gives rise to a graded group $h_*(X)$
(respectively, $h^{-\ast}(X)$) for every space $X$. The following
definition generalizes Shchepin's concept of total cohomology (see \cite{S}).


\begin{Def} \label{XX5.1} Let $G$ be an Abelian group.
For any pointed compact space $X$ its {\bf graded total cohomology
group} $\cal H^{-\ast}(X;G)$
is the direct sum of graded groups $H^{-\ast}(X,A;G)$
with $A$ ranging over all pointed closed subsets of $X$.
\end{Def} 

\begin{Thm} \label{XX5.2} Suppose $X$, $Y$ are pointed compact
spaces and $G$ is an Abelian group.
\par 1. $\cal H^{-\ast}(X)\wedge G\equiv \cal H^{-\ast}(X;G)$.
\par{2.} $\dim(\cal H^{-\ast}(X)\wedge \cal H^{-\ast}(Y))=\dim(\cal
H^{-\ast}(X\wedge Y))$.
\end{Thm} 
\begin{pf} 1) follows from the Universal coefficient formula
for cohomology.
\par
Recall that 4.11 in \cite{Dy$_6$} implies
 $$H^{-n}(X\wedge Y)\equiv \bigoplus\limits_{i}H^{-i}(X;H^{-(n-i)}(Y)).$$
Using that one can view $\cal H^{-\ast}(X))\wedge \cal H^{-\ast}(Y)$
 as a direct summand of $\cal H^{-\ast}(X\wedge Y)$
as follows:
$\cal H^{-\ast}(X)=\bigoplus\limits_{A\subset X}\cal H^{-\ast}(X/A)$,
$\cal H^{-\ast}(Y)=\bigoplus\limits_{B\subset Y}\cal H^{-\ast}(Y/B)$,
and $$(\cal H^{-\ast}(X)\wedge  \cal H^{-\ast}(Y))(n)\equiv $$
$$\bigoplus\limits_{A\subset X,B\subset Y}\bigoplus\limits_{i}H^{-i}(X/A;H^{-(n-i)}(Y/B))\equiv $$
$$\bigoplus\limits_{A\subset X,B\subset Y}H^{-n}(X\wedge Y/(X\wedge B\cup
A\wedge Y)).$$
\par Thus, 
$\dim(\cal H^{-\ast}(X)\wedge \cal H^{-\ast}(Y))\ge\dim(\cal H^{-\ast}(X\wedge Y))$
holds.
\par Suppose $\dim_G(\cal H^{-\ast}(X)\wedge \cal H^{-\ast}(Y))\ge n$
for some Abelian group $G$. In view of the above we get that
$H^i(X\times Y/(X\times B\cup A\times Y);G)=0$ for all $i\ge -n$,
 all closed subsets $A$ of $X$, and all closed subsets $B$ of $Y$.
Let $Z=X\times Y$.
Consider  
$${ \cal V}=\{   V\subset     Z \mid V \text{ is open and }
 H^i(Z,Z-V;G)=0 \text{ for all }i\ge  -n\}.$$
 Our goal is to show that $\cal
V$ contains all open sets in $Z$ which would complete the proof of 2).
Notice that $U\times V\in\cal V$ if $U$ is open in $X$ and $V$ is open in $Y$.
Indeed, $Z/(Z-U\times V)=X\times Y/((X-U)\times Y\cup X\times (Y-V))$.
The family $\{U\times V\mid U\text{ is open in }X, V\text{ is open in }Y\}$
is denoted by $\cal U$ and it is a subset of $\cal V$.
\par
Notice that if $V_1\subset     V_2\subset    \ldots$ is an increasing
sequence of elements in $\cal V$, then $\bigcup\limits_j V_j\in  \cal V$. It
is so as $H^i(Z,Z-\bigcup\limits_j V_j;G)=\text{dirlim} H^i(Z,Z-V_j;G)=0$
 for all $i\ge -n$.
Let us show that $V,W,V\cap W\in  \cal V$ implies $V\cup    W\in
{ \cal V}$. Let $A=Z\cup    C(Z-V)$ and $B=Z\cup    C(Z-W)$ be subsets of
the cone $C(Z)$ over $Z$. Then, $H^i(A;G)=H^i(Z,Z-V;G)=0$,
$H^i(B;G)=H^i(Z,Z-W;G)=0$ and $H^i(A\cup    B;G)=H^i(Z,Z-V\cap W;G)=0$
for all $i\ge -n$. From the Mayer-Vietoris exact sequence
 $H^i(A\cup    B;G)\to  H^i(A;G)\oplus     H^i(B;G)\to  H^i(A\cap
B;G)\to  H^{i+1}(A\cup    B;G)\to \ldots$ we get  $H^i(A\cap
B;G)= H^i(Z,Z-(V\cup    W);G)=0$ for all $i\ge -n$. Thus, $V\cup
W\in  { \cal V}$. Now, it is easy to show by induction on $m\ge  1$
that the union of $m$ elements from ${ \cal U}$ belongs to ${ \cal V}$.
Since every open set in $Z$ can be expressed as a union of an increasing
sequence of open sets, each of which is a finite union of elements of
$\cal U$, the proof of 2) is complete.
 \end{pf}

Recall that, for an unpointed compact space $X$, its
{\bf cohomological dimension} $\dim_G(X)$ can be defined in
two equivalent ways (see \cite{K}):
\par 1. As the smallest integer $n$
such that $H^k(X,A;G)=0$ for all $k\ge n+1$
and all closed subsets $A$ of $X$.
\par 2. As the smallest integer $n$ such that
$K(G,n)\in AE(X)$.
\par We will define $\dim_G(X)$ for pointed compact spaces
using the concept of homological dimension.

\begin{Def} \label{XX5.3} Suppose $X$ is a pointed compact
space. Given an Abelian group $G$
define {\bf the cohomological dimension}
$\dim_G(X)$ (also denoted by $d_X(G)$) as $-\dim_G(\cal H^{-\ast}(X))=
-\dim_\Z(\cal H^{-\ast}(X;G))$.
\end{Def} 

Notice that $\dim_0(X)=-\infty$ for any pointed compact space $X$.
Also, $\dim_G(X)=-\infty$ for any $G$ if $X$ is a pointed point.
However, in the remaining cases, our concept of the cohomological
dimension coincides with the one for the corresponding unpointed space.

\begin{Prop} \label{XXeh.2} Let $G$ be a non-trivial Abelian group.
If $X$ is a compact space containing at least two points,
then $\dim_G(X,x_0)=\dim_G(X)$ for all $x_0\in X$.
\end{Prop} 
\begin{pf} Notice that, for $n\ge 1$,
the groups $H^n(X,A;G)$ and $H^n(X,A\cup \{x_0\};G)$
are isomorphic. Consequently, $\dim_G(X)\leq n-1$
is equivalent to $\dim_G(X,x_0)\leq n-1$.
Hence, the only case where
$\dim_G(X,x_0)\ne \dim_G(X)$ may occur
is $\dim_G(X)=0$. In that case $X$ is totally disconnected,
so $H^0(X,x_0;G)\ne 0$ as $X$ contains at least two points.
Thus, $\dim_G(X,x_0)=0$ in that case as well.
\end{pf}

We generalize the concept of cohomological
dimension with respect to a group to that
of cohomological dimension with respect to a graded group
as it simplifies calculations.

\begin{Def} \label{XXeh.3} Suppose $X$ is a pointed compact
space. Given a graded group $A$
define {\bf the cohomological dimension}
$\dim_A(X)$ as $-\dim_\Z(\cal H^{-\ast}(X)\wedge A)$.
\end{Def}

\begin{Prop} \label{XXeh.4} Suppose $X$ is a pointed compact
space and $A$ is a graded group.
The following conditions are equivalent:
\par{1.} $\dim_A(X)\leq k$.
\par{2.} $\dim_{A(n)}(X)\leq n+k$ for each $n\in\Z$.
\end{Prop} 
\begin{pf} Let $B=\cal H^{-\ast}(X)$.
Since $A=\bigoplus\limits_{n\in\Z}\Sigma^n(A(n))$,
$$\dim_A(X)=-\dim_\Z(\bigoplus\limits_{n\in\Z}\Sigma^n(A(n)\wedge B)=$$
$$-\inf_{n\in\Z}\{\dim_\Z(\Sigma^n(A(n)\wedge B))\}=
-\inf_{n\in\Z}\{n-\dim_{A(n)}(X)\}.$$
Thus, $\dim_A(X)\leq k$ iff $n-\dim_{A(n)}(X)\ge -k$
for each $n$ which is equivalent to 2).
\end{pf}

\begin{Cor} \label{XX5.5} Suppose $X$ is a pointed compact
space and $K$ is a pointed CW complex.
The following conditions are equivalent:
\par{1.} $SP(K)\in AE(X)$.
\par{2.} $\cal H^{-\ast}(X)\wedge H_\ast(K)$ is non-negative.
\end{Cor} 
\begin{pf}  By \ref{XXX2.10b} $SP(K)\in AE(X)$ is equivalent to
$\dim_{H_n(K)}(X)\leq n$ for all $n\in\Z$.
By \ref{XXeh.4} that is equivalent to $\dim_{H_\ast(K)}(X)\leq 0$
which is another way of saying that 
$\cal H^{-\ast}(X)\wedge H_\ast(K)$ is non-negative in view of definition  \ref{XXeh.3}.
 \end{pf}

\begin{Cor} \label{XX5.7} Suppose $X$, $Y$ are pointed compact
spaces and $K$, $L$ are pointed CW complexes. If
$SP(K)\in AE(X)$ and $SP(L)\in AE(Y)$, then $SP(K\wedge L)\in AE(X\wedge Y)$.
\end{Cor} 

\begin{pf}  Suppose $Z$ is a compact space and
$M$ is a CW complex. \ref{XX5.5} says that $SP(M)\in AE(Z)$ if and only if
$\cal H^{-\ast}(Z)\wedge H_\ast(M)$ is non-negative. Thus, both
$\cal H^{-\ast}(X)\wedge H_\ast(K)$ and $\cal H^{-\ast}(Y)\wedge H_\ast(L)$ are non-negative.
Obviously, their smash product is non-negative
which means (see \ref{XX5.2}) that
$\cal H^{-\ast}(X\wedge Y)\wedge H_\ast(K\wedge L)$ is non-negative,
i.e., $SP(K\wedge L)\in AE(X\wedge Y)$ by \ref{XX5.5}.
 \end{pf}

\begin{Rem} \label{XXeh.5} Compare the simplicity
of the proof of \ref{XX5.7} with that of Theorem 5.4 of \cite{D-D$_2$}
($K\in AE(X)$ and $L\in AE(Y)$ imply
$SP(K\wedge L)\in AE(X\wedge Y)$ provided $X$ is metrizable,
$Y$ is metrizable and $\sigma$-compact, $K$ and $L$ are CW complexes).
\end{Rem}

\section{Bockstein groups and Bockstein basis}


\begin{Def} \label{XX4.11} Let $\cal{DGG}$ be the {\bf dimension algebra
of graded
groups}. It consists of
equivalence classes of graded groups: $A\sim B$ means $\dim(A)=\dim(B)$.
The direct product serves as the sum and the smash product
serves as the product in $\cal{DGG}$.

\end{Def}

Notice that $\Z$ serves as the unit of $\cal{DGG}$ ($A\wedge\Z=A$
for all graded groups $A$).
In this section we will:
\par{1.} show that all objects of $\cal{DGG}$ form a set,
\par{2.} find its basis.
\par

\begin{Def} \label{XX4.13} Let $\B$ (the {\bf Bockstein groups}) be the set of
groups consisting of the rationals $\Q$, the cyclic groups $\Z/\p$
for all primes $\p>0$, the quasi-cyclic groups $\Z/\p^{\infty}$
for all primes $\p>0$, and $\Z_{(\p)}$ ($\Z$ localized at $\p$)
for all primes $\p>0$.
\end{Def} 

\begin{Prop} \label{XXbg.0} Let $G$ be an Abelian group. If $0\to A\to B\to C\to 0$
is short exact sequence of graded groups, then
one has a long exact sequence of graded groups
$$\to\Sigma^{-1}(C)\wedge G\to A\wedge G\to B\wedge G\to C\wedge G\to \Sigma(A)\wedge G\to 
\Sigma(B)\wedge G\to \Sigma(C)\wedge G\to$$
If $C$ is torsion-free, then one has a short exact sequence
$$0\to A\wedge G\to B\wedge G\to C\wedge G\to 0.$$
\end{Prop} 
\begin{pf} By Corollary 9 on p.224 in \cite{Sp} one has an exact
sequence $0\to A\ast G\to B\ast G\to C\ast G\to A\otimes G\to B\otimes B\to C\otimes
G\to 0$. If $C$ is torsion-free, then it implies exactness
of $0\to A\wedge G\to B\wedge G\to C\wedge G\to 0.$
For arbitrary $C$ one splices $0\to A\ast G\to B\ast G\to C\ast G\to A\otimes G\to B\otimes B\to C\otimes
G\to 0$ to obtain 
$\to\Sigma^{-1}(C)\wedge G\to A\wedge G\to B\wedge G\to C\wedge G\to \Sigma(A)\wedge G\to 
\Sigma(B)\wedge G\to \Sigma(C)\wedge G\to$.
\end{pf}

\begin{Cor} \label{XXbg.1} If $0\to A\to B\to C\to 0$
is an exact sequence of graded groups then
\par 1. $A\oplus C\leq B$, $B\oplus \Sigma^{-1}(C)\leq A$, and $\Sigma(A)\oplus B\leq C$.
\par 2. If $C$ is torsion-free, then $\dim(B)=\dim(A\oplus C)$.
\par 3. If $\dim(A)=\dim(C)$ and $\dim_\Z(A)  > -\infty$, then $\dim(B)=\dim(A)$.
\end{Cor} 
\begin{pf} 1) and 2) follow from \ref{XXbg.0}.
\par Suppose $\dim(A)=\dim(C)$. Therefore $A\leq B$ by 1) and if $A < B$,
then there is $k$ and an Abelian group $G$ such that 
$\dim_G(B)\ge k+1$ but $\dim_G(A)=k$.
That however means $(\Sigma(A)\wedge G)(k)=0$
implying $(C\wedge G)(k)=0$ and contradicting $\dim_G(C)=\dim_G(A)=k$.
\end{pf}

Given any graded group $A$, $\Tor(A)$ is defined via $\Tor(A)(n)=\Tor(A(n))$
for each $n$. Also, $A/\Tor(A)$ is defined via $A/\Tor(A)(n)=A(n)/\Tor(A(n))$
for each $n$. Since $0\to \Tor(A)\to A\to A/\Tor(A)\to 0$
is exact and $A/\Tor(A)$ is torsion-free, \ref{XXbg.1}
implies the following.

\begin{Cor} \label{XX4.12} $\dim(A)=\dim((A/\Tor(A))\oplus \Tor(A))$
for any graded group $A$.
\end{Cor}

\begin{Cor} \label{XXbg.2} Let $\p$ be a prime.
\par 1. $\dim(\Z/\p^k)=\dim(\Z/\p)$ for all $k\ge 1$.
\par 2. $\Z/\p \leq H\leq \Sigma(\Z/\p)$ for all $\p$-groups $H$.
\par 3. $H\leq \Z/\p$ if and only if $H$ is not divisible by $\p$.
\par 4. $\dim(H)=\dim(\Z/\p)$ if $H$ is a $\p$-group
not divisible by $\p$.
\par 5. $\dim(H)=\dim(\Z/\p^\infty)$ if $H\ne 0$ is a $\p$-group
divisible by $\p$.
\par 6. $Z_{(\p)} < Z/\p < \Z/\p^\infty < \Sigma(\Z/\p)$.
\end{Cor} 
\begin{pf} Notice that, for every $k\ge 1$, there is an
exact sequence $0\to \Z/\p^k\to \Z/\p^{k+1}\to \Z/\p\to 1$.
Therefore, 1) follows from Part 3 of \ref{XXbg.1}.
\par By 1) all finite $\p$-groups $G\ne 0$ satisfy
$\dim(G)=\dim(\Z/\p)$. Since any $\p$-group is the direct limit
of its finite subgroups, 2) follows.
\par 3. If $\p\cdot H\ne H$, then $H\otimes(\Z/\p)=H/\p\cdot H$ is a direct sum of copies of $\Z/\p$.
Therefore $H\leq H\wedge(\Z/\p)\leq H\otimes(\Z/\p)\leq \Z/\p$.
If $H\leq \Z/\p$, then $H\wedge \Z/\p\leq \Z/\p\wedge \Z/\p\leq \Z/\p$
and $H/\p\cdot H=H\otimes \Z/\p=(H\wedge \Z/\p)(0)$ cannot be 0. 
\par 4) follows from 2) and 3).
\par 5. If $H\ne 0$ is a $\p$-divisible $\p$-group,
it is a direct sum of copies of $\Z/\p^\infty$.
\par 6. $Z_{(\p)} \leq Z/\p \leq \Z/\p^\infty \leq \Sigma(\Z/\p)$
follows from 2) and 3).
$\dim(Z_{(\p)})\ne  \dim(Z/\p)$ as $Z_{(\p)}\otimes\Z/\p^\infty\ne 0$
and  $\Z/\p\otimes\Z/\p^\infty=0$.
$\dim(\Z/\p^\infty)\ne  \dim(Z/\p)$ as $Z/\p\otimes\Z/\p\ne 0$
and  $\Z/\p\otimes\Z/\p^\infty=0$.
$\dim(\Z/\p^\infty)\ne\dim(\Sigma(\Z/\p))$
as $\dim_\Z(\Z/\p^\infty)=0\ne\dim_\Z(\Sigma(\Z/\p))=1$.
\end{pf}

\begin{Cor} \label{XXbg.3} If $G$ is a torsion-free group, then
$\dim(G)=\dim(\Q\oplus\bigoplus\limits_{\p\in l}\Z_{(\p)})$,
where $l=\{\p \mid \p\cdot G\ne G\}$.
\end{Cor} 
\begin{pf} It suffices to show
that, for any torsion-free group $F$, the condition $F\otimes H=0$
is equivalent to $H=\Tor(H)$ and $\p-\Tor(H)=0$ for any $\p$ not dividing $F$.
Clearly, if $\q\cdot F=F$, then $F\otimes H=0$ for any $\q$-group $H$.
If $\p\cdot F\ne F$ and $H\ne 0$ is a $\p$-group,
then, in view of $F\leq \Z/\p$ (see \ref{XXbg.2}) and $H\leq \Sigma(\Z/\p)$ (see \ref{XXbg.2}),
we get $F\otimes H=F\wedge H\leq (\Z/\p)\wedge \Sigma(\Z/\p)\leq \Sigma^2(\Z/\p)\ne 0$
which implies that $F\otimes H\ne 0$.
\end{pf}

\begin{Def} \label{XX4.14} Given an Abelian group $G$ define
its {\bf Bockstein basis} $\sigma(G)$ as the set of all Bockstein groups
$H$ such that $\dim(G)\leq \dim(H)$.

\end{Def} 

\begin{Cor} \label{XX4.16} $\dim(G)=\dim(\bigoplus\limits_{H\in\sigma(G)}H)$
for any Abelian group $G$.

\end{Cor}

\begin{pf}  Since $\dim(G)=\dim((G/\Tor(G))\oplus \Tor(G))$,
\ref{XXbg.2} and \ref{XXbg.3} imply that
there is a subset $S$ of $\sigma(G)$ such that 
$\dim(G)=\dim(\bigoplus\limits_{H\in S}H)$.
However, $G\leq H$ implies $\dim(G)=\dim(G\oplus H)$,
so \ref{XX4.16} follows.
 \end{pf}

\begin{Prop} \label{XX4.17} The following conditions are equivalent
for any graded Abelian groups $A$ and $B$:
\par{1.} $\dim(A)\leq \dim(B)$.
\par{2.} $\dim_G(A)\leq \dim_G(B)$ for all $G\in \B$.

\end{Prop}

\begin{pf}  Clearly, 1)$\implies$ 2).
If $\dim_G(A)\leq \dim_G(B)$ for all $G\in \B$,
then $\dim_G(A)\leq \dim_G(B)$ for all $G$ which are direct sums of Bockstein
groups. By \ref{XX4.16}, 2)$\implies$ 1).
 \end{pf}

The next result shows that our definition of the Bockstein basis coincides
with the one in \cite{Dy$_6$} and the only difference with the definitions in
\cite{K} or \cite{D$_8$} is that the case of $\Z/\p^\infty$ is treated
differently.

\begin{Prop} \label{XX4.15} Let $G$ be an Abelian group.
\par{1.} $\Q\in\sigma(G)$ if and only if $G\ne\Tor(G)$.
\par{2.} $\Z_{(\p)}\in\sigma(G)$ if and only if $G/\Tor(G)$
is not divisible by $\p$.
\par{3.} $\Z/\p\in\sigma(G)$ if and only if $G$
is not divisible by $\p$.
\par{4.} $\Z/\p^\infty\in\sigma(G)$ if and only if either
$\Z_{(\p)}\in\sigma(G)$
or $\p-\Tor(G)\ne 0$.

\end{Prop}

\begin{pf}  1. If $\dim(G)\leq\dim(\Q)$,
then $G\otimes \Q=G\wedge\Q\leq \Q\otimes\Q=\Q$,
so $G\otimes \Q\ne 0$ and $G$ cannot be a torsion group.
If $G\ne\Tor(G)$, then (see \ref{XX4.12}) $\dim(G)=\dim(\Tor(G)\oplus G/\Tor(G))\leq
\dim(G/\Tor(G))\leq\dim(\Q)$ by \ref{XXbg.3}.
\par 2. Suppose $\dim(G)\leq\dim(\Z_{(\p)})$
and $G/\Tor(G)$ is divisible by $\p$. 
Therefore $0=G\otimes \Z/\p^\infty$ and since
$G\wedge\Z/\p^\infty\leq \Z_{(\p)}\wedge\Z/\p^\infty=\Z/\p^\infty$,
we arrive at a contradiction.
If $G/\Tor(G)$ is not divisible by $\p$, then (see \ref{XXbg.3})
$G/\Tor(G)\leq \Z_{(\p)}$. Since $G\leq G/\Tor(G)$,
$\Z_{(\p)}\in\sigma(G)$.
\par 3. Follows from Part 3 of \ref{XXbg.2}.
\par 4. Since $Z_{(\p)} < \Z/\p^\infty$,
$Z_{(\p)}\in\sigma(G)$ implies $\Z/\p^\infty\in\sigma(G)$.
If $\p-\Tor(G)\ne 0$, then $G\leq \Tor(G)\leq \p-\Tor(G)\leq \Z/\p^\infty$
by \ref{XXbg.2}.
Suppose $\dim(G)\leq\dim(\Z/\p^\infty)$, $G/\Tor(G)$ is divisible
and $\p-\Tor(G)=0$. Hence $\p\cdot G=G$ and $G\ast \Z/\p=0$.
Consequently $0=G\wedge \Z/\p\leq \Z/\p^\infty\wedge\Z/\p\ne 0$,
a contradiction.
 \end{pf}

\begin{Def} \label{XX4.18} Given a function $\alpha:\B\to\Z\cup\{\pm\infty\}$
define $GG(\alpha)$ as $\bigoplus\limits_{H\in\B}\Sigma^{\alpha(H)}(H)$.

\end{Def}


\begin{Prop} \label{XX4.19} For every Abelian group $A$
there is a function $\alpha:\B\to\Z\cup\{\pm\infty\}$
such that $\dim(A)=\dim(GG(\alpha))$.

\end{Prop}

\begin{pf}  Given a graded group $B$ we define a new graded group $B'$ via
$B'(n)=\bigoplus\limits_{i\leq n}B(i)$ and notice that $\dim(B)=\dim(B')$.
\par
Given a Bockstein group $H$ there are three possibilities:
\par{a.} $H\in\sigma(A'(n))$ for each $n$. Put $\alpha(H)=-\infty$.
\par{b.} $H\notin\sigma(A'(n))$ for each $n$. Put $\alpha(H)=\infty$.
\par{c.} $H\in\sigma(A'(n))$ and $H\notin\sigma(A'(n-1))$ for some $n\in\Z$. Put $\alpha(H)=n$.
\par Let $B=GG(\alpha)$ and notice that
$B'(n)=\bigoplus\limits_{H\in\sigma(A'(n))}H$ which implies
$\dim(B')=\dim(A')$ by \ref{XX4.16}.
 \end{pf}

The above results mean that every graded group has the same dimension
as a graded group whose terms are direct sums of Bockstein groups.
Therefore, $\cal{DGG}$
forms a set. In the future we will switch to graded groups with terms
being direct sums of Bockstein groups.

\section{Hurewicz theorems in extension theory}

This section deals with analogs of Hurewicz Theorem in extension theory.
The most important case is that of simply connected CW complexes.
However, our proofs work for a wider class of spaces, namely nilpotent
CW complexes, so that is the setup we chose for this section.
\par
Let us start with a result that is closest in spirit to the classical
Hurewicz Theorem.
\begin{Lem} \label{XXhur.1} 
Let $K$ be a nilpotent pointed CW complex such that $\pi_\ast(K)$ is Abelian.
If 
 $G$ is a Bockstein group and $\dim_G(H_\ast(K))\ge k$,
then 
$(\pi_\ast(K)\wedge G)(k)$ is isomorphic to $(H_\ast(K)\wedge G)(k)$.
\end{Lem} 
\begin{pf} If $k\leq 1$, then \ref{XXhur.1} is obvious
in view of $\pi_1(K)=H_1(K)$. Therefore only $k\ge 2$ is of
interest.
\par
Case 1: Consider $G=\Z_{(l)}$, where $l$ is a set of primes.
Recall that $\Z_{(l)}$ is the subring of rationals $\Q$
consisting of ratios $m/n$, where $n$ is not divisible
by all $\p\in l$. In particular, $\Z_{(\emptyset)}=\Q$.
Use a localizing map $f:K\to K_G$ of pointed CW complexes
such that $\pi_\ast(f)$ is $\pi_\ast(K)\to \pi_\ast(K)\otimes G$
and $H_\ast(f)$ is $H_\ast(K)\to H_\ast(K)\otimes G$
(see  \cite{BK} or \cite{HMR}, also \cite{Su} for $K$ simply connected).
Now \ref{XXhur.1} follows from the classical Hurewicz Theorem.
\par Case 2: $G=\Z/\p^\infty$ for some prime $\p$.
The exact sequence $0\to \Z\to \Z[1/\p]\to G\to 0$
yields $0\to H\ast G\to H\to H\otimes \Z[1/\p]\to H\otimes G\to 0$
for any Abelian group $H$.
Consider the localizing map $f:K\to L=K_{1/\p}$
such that $\pi_\ast(f)$ is $\pi_\ast(K)\to \pi_\ast(K)\otimes \Z[1/\p]$
and $H_\ast(f)$ is $H_\ast(K)\to H_\ast(K)\otimes \Z[1/\p]$.
We may assume that $K$ is a subcomplex of $L$
and $f$ is the inclusion map.
The exact sequence
$$\pi_m(K)\to \pi_m(L)\to \pi_m(L,K)\to \pi_{m-1}(K)\to \pi_{m-1}(L)$$
implies existence of an exact sequence
$$0\to \pi_{m}(K)\otimes \Z/\p^\infty\to \pi_{m}(L,K)\to \pi_{m-1}(K)\ast \Z/\p^\infty\to 0$$
as the cokernel of $\pi_m(K)\to \pi_m(L)$ is $\pi_{m}(K)\otimes \Z/\p^\infty$
and the kernel of $\pi_{m-1}(K)\to \pi_{m-1}(L)$ is $\pi_{m-1}(K)\ast \Z/\p^\infty$.
Since $\pi_{m}(K)\otimes \Z/\p^\infty$ is a divisible group, it is a direct summand
of $\pi_m(L,K)$ and the above sequence splits. In particular,
$\pi_{m}(L,K)\equiv (\pi_\ast(K)\wedge \Z/\p^\infty)(m)$ for all $m$.
Similarly, $H_{m}(L,K)\equiv (H_\ast(K)\wedge \Z/\p^\infty)(m)$ for all $m$.
Therefore $\pi_m(L,K)=0$ for $m\leq k-1$
and the Hurewicz homomorphism $\pi_m(L,K)\to H_m(L,K)$
is an isomorphism for $m=k$ and an epimorphism
for $m=k+1$.
\par Case 3:  $G=\Z/\p$ for some prime $\p$.
We will use mod $\p$ homotopy groups
$\pi_n(K;\Z/\p)$ constructed in \cite{N}.
In view of Proposition 1.4 on p.3 in \cite{N} 
one has an exact sequence
$$ 0\to \pi_k(K)\otimes \Z/\p\to \pi_k(K;\Z/\p)\to \pi_{k-1}(K)\ast \Z/\p\to 0$$
and the mod $\p$ Hurewicz Theorem in \cite{N} (see 3.8 on p.12)
implies that $\pi_k(K;\Z/\p)\equiv H_k(K;\Z/\p)$.
Therefore $\pi_k(K;\Z/\p)$ is a vector space over $\Z/\p$
and the above exact sequence splits
showing $(\pi_\ast(K)\wedge G)(k)$ being isomorphic to $(H_\ast(K)\wedge G)(k)$.
\end{pf}

The following result is the algebraic version
of the Hurewicz theorem in extension theory (for a more
geometric version see \ref{XXeh.9}).
In the case of Abelian CW complexes it has been proved
by Shchepin \cite{S} under a different form
as a stronger result than \ref{XXeh.9}. In our exposition
both versions are equivalent in view of \ref{XXeh.4}.

\begin{Cor} \label{XXeh.8} 
If $K$ is a nilpotent pointed CW complex such that $\pi_\ast(K)$ is Abelian,
then $$\dim(H_\ast(K))=\dim(\pi_\ast(K)).$$
\end{Cor} 
\begin{pf}
In view of  \ref{XX4.17} it suffices to show
that $\dim_H(H_\ast(K))=\dim_H(\pi_\ast(K))$
for all Bockstein groups $H$ which follows from  \ref{XXhur.1}.
\end{pf}

\begin{Prop} \label{XXeh.6} 
Suppose $A$ and $B$ are two graded groups such that $A\leq B$ does not hold.
If $A$ is positive,
then there is a pointed compactum $X$ of finite dimension such that
$\dim_A(X)\leq 0 < \dim_B(X)$.
\end{Prop} 
\begin{pf} We may assume that $A(n)$ is countable for each $n$.
Since $A \leq B$ is false, there is an Abelian group $G$
such that $\dim_G(A) > \dim_G(B)=k$. 
Let $K$ be the wedge of $S^{k+1}$ and Moore spaces $M(A(n),n)$ for $n\leq k$.
Notice that $H_n(K;G)=0$ for $n\leq k$.
If $B(k)\oplus G\ne 0$, put $P=K(B(k),k+1)$.
If $B(k-1)\ast G\ne 0$, put $P=K(B(k-1),k)$.
We are going to use Theorem II of \cite{Dy$_4$}:
\par
Suppose $G$ is an abelian group, $m>0$ and $K$
is a countable connected CW complex. Then, the following 
conditions are equivalent:
\par{1.} For any CW complex $P$ and any $a\in H_m(P;G)-\{0\}$
there is a compactum $X$ and a map $\pi :X\to P$
such that 
$$a\in \im (\check H_m(X;G)\to \check H_m(P;G))$$
and $K$ is an absolute extensor of $X$.
\par{2.} $\Tilde H_k(K;G)=0$ for all $k<m$.
\par Since $H_{k+1}(P;G)\ne 0$,
there is a non-trivial map $X\to P$
such that $K\in AE(X)$.
In particular, $\dim(X)\leq k+1$ which implies
$\dim_{A(i)}(X)\leq i$ for all $i\ge k+1$.
Since $M(A(i),i)\in AE(X)$ for $i\leq k$,
we get $K(A(i),i)\in AE(X)$ (see \cite{D$_5$}).
By  \ref{XXeh.4} one gets $\dim_A(X)\leq 0$.
However, since $f:X\to P$ is non-trivial,
either $\dim_{B(k)}(X)>k$ or $\dim_{B(k-1)}(X)>k-1$.
Therefore, $\dim_B(X)\leq 0$ does not hold (see \ref{XXeh.4}).
 \end{pf}

\begin{Cor} \label{XXeh.7} 
The following conditions are equivalent for any positive
graded groups
$A$ and $B$:
\par 1. $\dim(A)\leq \dim(B)$.
\par 2. $\dim_A(X)\ge \dim_B(X)$ for all
pointed compact spaces $X$.
\par 3. $\dim_A(X)\ge \dim_B(X)$ for all
finitely dimensional pointed compacta $X$.
\par 4. $\dim_A(X)\leq 0$ implies $\dim_B(X)\leq 0$ for every finite-dimensional
pointed compactum $X$.
\end{Cor} 
\begin{pf} 1)$\implies$2). Put $C=\cal H^{-\ast}(X)$.
Since $A\leq B$, $A\wedge C\leq B\wedge C$
and $\dim_B(X)=-\dim_\Z(B\wedge C)\leq -\dim_\Z(A\wedge C)=\dim_A(X)$.
\par Both 2)$\implies$3) and 3)$\implies$4) are obvious.
\par 4)$\implies$1) follows from \ref{XXeh.6}. 
 \end{pf}

The following was proved in \cite{D$_5$} in the simply connected case
and in \cite{CD} in the nilpotent case.

\begin{Thm} \label{XXeh.9} 
If $K$ is a nilpotent pointed CW complex
and $X$ is a compact space of finite dimension,
then the following conditions are equivalent:
\par 1. $K\in AE(X)$.
\par 2. $SP(K)\in AE(X)$.
\par 3. $\dim_{H_n(K)}(X)\leq n$ for all $n\ge 1$.
\par 4. $\dim_{\pi_n(K)}(X)\leq n$ for all $n\ge 1$.
\end{Thm} 
\begin{pf} 1)$\implies$2) is the same as \ref{XXX2.10a}.
\par 2)$\implies$3) is the same as \ref{XXX2.10b}.
\par 3)$\implies$4) follows from \ref{XXeh.8}, \ref{XXeh.7}, and \ref{XXeh.4}.
\par 4)$\implies$1) follows from \ref{XXXintro.1}.
 \end{pf}

\section{The dual of graded groups}

In this section we will define the dual $A^*$ for each $A$.
 $A^*$ is the algebraic manifestation of the duality
between
compact spaces $X$ and their extension dimension $\ExD(X)$ (which is
represented by
a CW complex).


\begin{Def} \label{XX4.20} Given a graded group $A$ define its {\bf dual}
$A^*$ as
\par\centerline{$\bigoplus\{\Sigma^k(H)\mid H\in\B\text{ and }
k\ge-\dim_H(A)\}.$}

\end{Def}


\begin{Prop} \label{XX4.21} $A^*$
is the minimum of $\{B\in\cal{DGG}\mid \dim_\Z(A\wedge B)\ge 0\}$.

\end{Prop}

\begin{pf}  Suppose $\dim_\Z(A\wedge B)\ge 0$ and
$B=\bigoplus\limits_{H\in\B}\Sigma^{\alpha(H)}(H)$.
Since $\dim(B)\leq \dim(\Sigma^{\alpha(H)}(H))$ for each $H\in\B$,
$\dim_\Z(A\wedge \Sigma^{\alpha(H)}(H))\ge 0$ and
$\alpha(H)+\dim_\Z(A\wedge H)=\alpha(H)+\dim_H(A)\ge 0$.
Thus, $\alpha(H)\ge -\dim_H(A)$ and $\dim(A^*)\leq
\dim(\Sigma^{\alpha(H)}(H))$
for each $H\in\B$ which means
$\dim(A^*)\leq\dim(\bigoplus\limits_{H\in\B}\Sigma^{\alpha(H)}(H))=\dim(B)$.
 \end{pf}


\begin{Cor} \label{XX4.22} If $\dim(A)\leq \dim(B)$, then $\dim(A^*)\ge
\dim(B^*)$.

\end{Cor}

\begin{pf}  $0\leq \dim_\Z(A\wedge A^*)\leq \dim_\Z(B\wedge A^*)$
implies $\dim(A^*)\ge \dim(B^*)$ by the previous result.
 \end{pf}


\begin{Thm} \label{XX4.23} $\dim((A^*)^*)=\dim(A)$.

\end{Thm}

\begin{pf}  $0\leq \dim_\Z(A\wedge A^*)$
implies $\dim((A^*)^*)\leq \dim(A)$ by \ref{XX4.21}.
\par
We need to show $\dim(A)\leq \dim((A^*)^*)$. Let $B=(A^*)^*$.
It suffices to show that $\dim_H(A)\ge k$
implies $\dim_H((A^*)^*)\ge k$ for all $H\in\B$.
Since $\dim_H(A)\ge k$,
$H$ is a direct summand of $A^*(-k)$.
If $F$ is a direct summand of $B(s)$ for $s\leq k-1$,
then $F\otimes A^*(-k)=0$ and $F*A^*(-k)=0$ if $s<k-1$.
Therefore, $H\otimes B(s)=0$ for $s\leq k-1$,
and $H*B(s)=0$ if $s<k-1$. This proves $\dim_H(B)\ge k$.
 \end{pf}

Given a graded group $A$ we are presented with two problems:
\par{a.} to represent $A$ as $GG(\alpha)$
for some optimal $\alpha:\B\to \Z\cup\{\pm\infty\}$,
\par{b.} to compute dimensions $\dim_H(A)$ of $A$ for all $H\in\B$.
\par\noindent It turns out
these two problems are intertwined. The natural choice for $\alpha(H)$ in a)
would be to seek the smallest integer $n$ so that
$\dim(A)\leq\dim(\Sigma^n(H))$.
This leads to the following concepts.


\begin{Def} \label{XX4.24} Given a graded group $A$ the
function $d_A:\cal{DGG}\to \Z\cup\{\pm\infty\}$ is defined by
$d_A(X)=\dim_\Z(A\wedge X)$. The function $e_A:\cal{DGG}\to
\Z\cup\{\pm\infty\}$
is defined as the infimum of all $m$ so that $\dim(A)\leq \dim(\Sigma^m(X))$.

\end{Def}


\begin{Thm} \label{XX4.25} $e_A(X)=-\dim_\Z(A^*\wedge X)=-d_{A^*}(X)$
 for all graded groups $A$ and $X$.

\end{Thm}

\begin{pf}  Suppose $e_A(X)\leq m$, i.e. $\dim(A)\leq\dim(\Sigma^m(X))$.
This implies
$0\leq \dim_\Z(A\wedge A^*)\leq\dim_\Z(\Sigma^m(X))\wedge
A^*)=m+\dim_\Z(X\wedge A^*)$,
i.e. $-\dim_\Z(X\wedge A^*)\leq m$. Thus, $-\dim_\Z(X\wedge A^*)\leq e_A(X)$.
\par Suppose $-\dim_\Z(X\wedge A^*)\leq k$. This implies
$0\leq k+\dim_\Z(X\wedge A^*)=\dim_\Z(\Sigma^k(X)\wedge A^*)$.
Therefore $\dim((A^*)^*)\leq\dim(\Sigma^k(X))$ and $e_A(X)\leq k$ as
$\dim((A^*)^*)=\dim(A)$. This means $e_A(X)\leq-\dim_\Z(X\wedge A^*)$.
 \end{pf}


\begin{Def} \label{XX4.26} A function $\alpha:\cal{DGG}\to\Z\cup\{\pm\infty\}$
is called a {\bf dimension-like function} if the following conditions are
satisfied:
\par{1.} $\alpha(\Sigma^k(A))=k+\alpha(A)$.
\par{2.} $\dim(A)\leq\dim(B)$ implies $\alpha(A)\leq\alpha(B)$.
\par{3} $\alpha(\bigoplus\limits_{t\in T}A_t)=\inf\limits_{t\in
T}\alpha(A_t)$.
\par $\alpha$ is called an {\bf extension function} if
$-\alpha$ is a dimension-like function.

\end{Def}


\begin{Cor} \label{XX4.27} $d_A$ is a dimension-like function
and $e_A$ is an extension function for every graded group $A$.

\end{Cor}

\begin{pf}  $d_A$ is clearly a dimension-like function for all $A$
which implies that $e_A$
is an extension function for all $A$.
 \end{pf}

Let us show that extension/dimension-like
 functions are completely determined by their values on
Bockstein groups.


\begin{Prop} \label{XX4.28} Given two extension functions
$e_1,e_2:\cal{DGG}\to\Z\cup\{\pm\infty\}$ the following conditions are
equivalent:
\par{1.} $e_1\leq e_2$.
\par{2.} $e_1(H)\leq e_2(H)$ for each $H\in\B$.

\end{Prop}

\begin{pf}  Only 2)$\implies$ 1) is of interest.
Let $d_1=-e_1$ and $d_2=-e_2$ be the corresponding dimension-like functions.
We know that $d_2(H)\ge d_1(H)$ for each $H\in\B$ and we need to show
$d_2(A)\ge d_1(A)$ for each graded group $A$.
First assume $A=\Sigma^k(H)$ for some $H\in\B$.
By 1) of \ref{XX4.26}, $d_2(A)=k+d_2(H)\ge k+d_1(H)=d_1(A)$.
By \ref{XX4.19}, any $A$ is equivalent to $\bigoplus\limits_{t\in T}A_t$
so that $d_2(A_t)\ge d_1(A_t)$ for each $t\in T$.
Now, 4.26.3 says that $d_2(A)=\inf\limits_{t\in
T}d_2(A_t)\ge \inf\limits_{t\in
T}d_1(A_t)=d_1(A)$.
 \end{pf}

A natural question arises which functions $\alpha:\B\to\Z\cup\{\pm\infty\}$
give rise to extension/dimension-like functions.


\begin{Def} \label{XX4.30} $\alpha:\B\to\Z\cup\{\pm\infty\}$
is a Bockstein function if the following conditions are satisfied:
\par{1.} $\alpha(\Z/{\p^{\infty}})\le \alpha(\Z/\p)$,
\par{2.} $\alpha(\Z/\p)\le \alpha(\Z/{\p^{\infty}})+1$,
\par{3.} $\alpha(\Q)\le \alpha(\Z_{(\p)})$,
\par{4.} $\alpha(\Z/\p)\le \alpha(\Z_{(\p)})$,
\par{5.} $\alpha(\Z/{\p^{\infty}})\le
\max(\alpha(\Q),\alpha(\Z_{(\p)})-1)$,
\par{6.} $\alpha(\Z_{(\p)})\le
\max(\alpha(\Q),\alpha(\Z/{\p^{\infty}})+1)$.

\end{Def}


\begin{Prop} \label{XX4.31} If $e:\cal{DGG}\to\Z\cup\{\pm\infty\}$
is an extension function, then its restriction
 $e|\B:\B\to\Z\cup\{\pm\infty\}$ is a Bockstein function.

\end{Prop}

\begin{pf}  We need to check conditions 1)-6) of \ref{XX4.30}.
1), 3), and 4) follow from \ref{XX4.15}.
Notice that $\Z/\p$ is a subgroup of $\Z/{\p^{\infty}}$
and one has an exact sequence
$0\to \Z/\p\to \Z/{\p^{\infty}}\to \Z/{\p^{\infty}}\to 0$.
Thus, 2) follows from  \ref{XXbg.1}.
Notice that $\Z_{(\p)}$ is a subgroup of $\Q$
and one has an exact sequence
$0\to \Z_{(\p)}\to \Q\to \Z/{\p^{\infty}}\to 0$.
Thus, 5) and 6) follow from  \ref{XXbg.1}.
 \end{pf}


\begin{Prop} \label{XX4.32} The following conditions are equivalent
for all graded groups $A$ and $B$:
\par{1.} $\dim(A)\leq \dim(B)$.
\par{2.} $d_A\leq d_B$.
\par{3.} $e_A\leq e_B$.

\end{Prop}

\begin{pf}  1)$\implies$ 2) follows from \ref{XX4.27}.
1) is a special case of 2) as $d_A(G)=\dim_G(A)$.
Since $e_A=-d_{A^*}$, 3) is equivalent to $\dim(B^*)\leq\dim(A^*)$ which
is equivalent to 1) by \ref{XX4.22} and \ref{XX4.23}.
 \end{pf}

\section{Algebra of Bockstein functions}

Our goal is to show that $\cal{DGG}$ is isomorphic to an algebra composed
of Bockstein functions (an equivalent approch would be to use
dimension-like functions; Bockstein functions are traditionally used
in cohomological dimension theory). To do that we need to show that for
every Bockstein function
$\beta$ there is a graded group whose extension function equals $\beta$.


\begin{Def} \label{XX4.33}
Given a graded Abelian group $A$
its
{\bf Bockstein graded group} $BGG(A)$ is defined
as $GG(e_A)$.

\end{Def}


\begin{Prop} \label{XX4.34} If $\beta:\B\to\Z\cup\{\pm\infty\}$ is a
Bockstein function
and $A=GG(\beta)$, then the following conditions are satisfied:
\par{1.} $e_A=\beta$.
\par{2.} $d_A(H)=e_A(H)$ for $H=\Q,\Z/\p$.
\par{3.} If $\beta(\Z_{(p)})=\beta(\Z/\p^\infty)$,
then $d_A(H)=e_A(H)$ for $H=\Z_{(p)},\Z/\p^\infty$.
\par{4.} If $\beta(\Z_{(p)})>\beta(\Z/\p^\infty)$,
then $d_A(\Z_{(p)})=\min(\beta(\Q),\beta(\Z/\p^\infty))$
and $d_A(\Z/\p^\infty)=\beta(\Z/\p^\infty)+1$.

\end{Prop}

\begin{pf}  Clearly, $e_A\leq\beta$.
Suppose $H\in\B$ and $\beta(H)=m$.
\par
If $H=\Q,\Z/\p$, then $H\otimes F=0$ for every $F\in\B$
such that $\beta(F)<\beta(H)$ and $0=H*F$ for every $F\in\B$
such that $\beta(F)<\beta(H)-1$. Since $H\otimes H\ne 0$,
$d_A(H)=e_A(H)=m$ in this case.
\par Assume $\beta(\Z_{(p)})=\beta(\Z/\p^\infty)=k$
(which implies $\beta(\Q)=k$) and let $G=\Z_{(p)}$ or $G=\Z/\p^\infty$.
Notice that
 $G\otimes F=0=G*F=0$ for every $F\in\B$
such that $\beta(F)<k$
which implies $d_A(H)=k$ if $H=\Z_{(p)}$ or $H=\Z/\p^\infty$. Also,
assuming $G=\Z_{(p)}$,
$\dim_G(A)=k>k-1=\dim_G(\Sigma^{k-1}(H))$
if $H=\Z_{(p)}$ or $H=\Z/\p^\infty$.
Thus, $e_A(H)=\beta(H)$ if $H=\Z_{(p)}$ or $H=\Z/\p^\infty$.
\par Assume $\beta(\Z_{(p)})=\beta(\Z/\p^\infty)+1=k$ which implies
$\beta(\Q)\leq k$.
Let $F=\Z_{(p)}$ and $G=\Z/\p^\infty$. Notice that $\dim_G(A)=k$
and $\dim_F(A)=\min(\beta(\Q),\beta(\Z/\p^\infty))$. Since
 $\dim_G(\Sigma^{k-1}(\Z_{(p)}))=k-1$ and
$\dim_G(\Sigma^{k-2}(\Z/\p^\infty))=k-1$,
we get $e_A(\Z_{(p)})=\beta(\Z_{(p)})$
and $e_A(\Z/\p^\infty)=\beta(\Z/\p^\infty)$.
\par Assume $k=\beta(\Z_{(p)})>\beta(\Z/\p^\infty)+1=m$ which implies
$\beta(\Q)=k$.
Let $F=\Z_{(p)}$ and $G=\Z/\p^\infty$. Notice that $\dim_G(A)=m$,
$\dim_F(A)=m-1$,
and $\dim_\Q(A)=k$. Since
 $\dim_\Q(\Sigma^{k-1}(\Z_{(p)}))=k-1$ and
$\dim_G(\Sigma^{m-2}(\Z/\p^\infty))=m-1$,
we get $e_A(\Z_{(p)})=\beta(\Z_{(p)})$
and $e_A(\Z/\p^\infty)=\beta(\Z/\p^\infty)$.
 \end{pf}


\begin{Cor} \label{XX4.35} $\dim(A)=\dim(BGG(A))$ for any graded Abelian
group $A$.

\end{Cor}

\begin{pf}  Let $B=BGG(A)$. By \ref{XX4.34}, $e_A=e_B$ which implies
$\dim(A)=\dim(B)$ by \ref{XX4.32}.
 \end{pf}


\begin{Cor} \label{XX4.36} If $A$ is graded Abelian group
and $\p$ is a prime number, then the following conditions are equivalent:
\par{1.} $e_A(\Z/\p^\infty)=e_A(\Z_{(p)})$.
\par{2.} $d_A(\Z/\p^\infty)=d_A(\Z_{(p)})$.
\par{3.} $e_A(\Z/\p^\infty)=d_A(\Z/\p^\infty)$.
\par{4.} $e_A(\Z_{(p)})=d_A(\Z_{(p)})$.

\end{Cor}

\begin{pf}  Let $B=BGG(A)$. By \ref{XX4.35},
$\dim(A)=\dim(B)$ which implies \ref{XX4.36} in view of \ref{XX4.34}.
 \end{pf}


\begin{Def} \label{XX4.37} A Bockstein/dimension-like function $\beta$ is
$\p$-regular
(respectively, $\p$-singular) if $\beta(\Z/\p^\infty)=\beta(\Z_{(p)})$
(respectively, $\beta(\Z/\p^\infty)\ne\beta(\Z_{(p)})$).
\par
A graded group $A$ is $\p$-regular
(respectively, $\p$-singular) if $e_A$ is $\p$-regular
(respectively, $\p$-singular).

\end{Def}


\begin{Cor} \label{XX4.38} If $A$ is a graded group,
then the following conditions are satisfied:
\par{1.} $d_A(H)=e_A(H)$ for $H=\Q,\Z/\p$.
\par{2.} If $A$ is $\p$-regular,
then $d_A(H)=e_A(H)$ for $H=\Z_{(p)},\Z/\p^\infty$.
\par{3.} If $A$ is $\p$-singular,
then $d_A(\Z_{(p)})=\min(e_A(\Q),e_A(\Z/\p^\infty))$,
 $d_A(\Z/\p^\infty)=e_A(\Z/\p^\infty)+1$,
and $e_A(\Z_{(p)})=\max(d_A(\Q),d_A(\Z/\p^\infty))$.

\end{Cor}

\begin{pf}  Let $B=BGG(A)$. By \ref{XX4.35},
$\dim(A)=\dim(B)$ which implies \ref{XX4.38} in view of \ref{XX4.34}.
 \end{pf}

To be able to add Bockstein functions in a meaningful way,
we assume the convention $+\infty+(-\infty)=+\infty$
as $-\infty$ is understood to be the infimum of all integers
(see the proof below).


\begin{Cor} \label{XX4.39} If $A, B$ are graded groups and $C=A\wedge B$,
then the following conditions are satisfied:
\par{1.} $d_C(H)=e_C(H)=e_A(H)+e_B(H)$ for $H=\Q,\Z/\p$.
\par{2.} If $A$ or $B$ is $\p$-regular,
then $d_C(H)=d_A(H)+d_B(H)$ and $e_C(H)=e_A(H)+e_B(H)$ for
$H=\Z_{(p)},\Z/\p^\infty$.
\par{3.} If both $A$ and $B$ are $\p$-singular,
then $$e_C(\Z/{\p^{\infty}})=\min(e_A(\Z/{\p})+
e_B(\Z/{\p}),e_A(\Z/{\p^{\infty}})+e_B(\Z/{\p^{\infty}})+1),$$
$$d_C(\Z/{\p^{\infty}})=\min(d_A(\Z/{\p})+
d_B(\Z/{\p})+1,d_A(\Z/{\p^{\infty}})+d_B(\Z/{\p^{\infty}})),$$
$$e_C(\Z_{(\p)})
=\max(e_C(\Q),e_C(\Z/{\p^{\infty}})+1),$$
$$d_C(\Z_{(\p)})
=\min(d_C(\Q),d_C(\Z/{\p^{\infty}})-1).$$

\end{Cor}

\begin{pf}  By \ref{XX4.35}, we may assume $A=GG(\alpha)$ and $B=GG(\beta)$,
where $\alpha=e_A$ and $\beta=e_B$.
We will concentrate on computing $d_C$ as the equations for $e_C$
follow from \ref{XX4.38}. Notice that
$C=\bigoplus\limits_{F,G\in\B}\Sigma^{\alpha(F)+\beta(G)}(F\wedge G)$.
Therefore, $d_C(H)=\inf\{\alpha(F)+\beta(G)+\dim_H(F\wedge G)\mid F,G\in\B\}$.
\par
If $H=\Q,\Z/\p$, then the infimum of $\alpha(F)+\beta(G)+\dim_H(F\wedge G)$
is achieved for $F=G=H$ (see \ref{XX4.30}) and is equal to $\alpha(H)+\beta(H)=
d_A(H)+d_B(H)$.
\par Assume $A$ is $\p$-regular.
If $H=\Z_{(p)}$, then the infimum of $\alpha(F)+\beta(G)+\dim_H(F\wedge G)$
is achieved either for $F=G=\Q$ (and equals $d_A(\Q)+d_B(\Q)=
d_A(\Z_{(p)})+d_B(\Q)$)
or is achieved
for $F=\Z_{(p)}$, $G=\Z/\p^\infty$ (see \ref{XX4.30}) and is equal to
$\alpha(F)+\beta(G)
=d_A(\Z_{(p)})+e_B(\Z/\p^\infty)$. Thus,
$d_H(C)=d_A(\Z_{(p)})+\min(d_B(\Q),e_B(\Z/\p^\infty))=d_A(\Z_{(p)})+d_B(\Z_{(p)}
)$.
If $H=\Z/\p^\infty$, then the infimum of $\alpha(F)+\beta(G)+\dim_H(F\wedge G)$
is achieved either for $F=\Z_{(p)}$, $G=\Z/\p^\infty$ (see \ref{XX4.30}) and is
equal to
$\alpha(F)+\beta(G)+1$, or for
$F=\Z_{(p)}$, $G=\Z_{(p)}$ (see \ref{XX4.30}) and is equal to
$\alpha(F)+\beta(G)$.
Thus,
$d_H(C)=d_A(\Z_{(p)})+\min(e_B(\Z_{(p)}),e_B(\Z/\p^\infty)+1)=d_A(\Z/\p^\infty)+
d_B(\Z/\p^\infty)$.
\par Assume both $A$ and $B$ are $\p$-singular.
If $H=\Z/\p^\infty$, then the infimum of $\alpha(F)+\beta(G)+\dim_H(F\wedge G)$
is achieved either for $F=\Z/\p$, $G=\Z/\p$ (see \ref{XX4.30}) and is equal to
$\alpha(F)+\beta(G)+1=d_A(\Z/\p)+d_B(\Z/\p)+1$, or for
$F=\Z/\p^\infty$, $G=\Z/\p^\infty$ (see \ref{XX4.30}) and is equal to
$\alpha(F)+\beta(G)+2=d_A(H)+d_B(H)$.
Thus, $d_C(\Z/{\p^{\infty}})=\min(d_A(\Z/{\p})+
d_B(\Z/{\p})+1,d_A(\Z/{\p^{\infty}})+d_B(\Z/{\p^{\infty}}))$.
If $H=\Z_{(p)}$, then the infimum of $\alpha(F)+\beta(G)+\dim_H(F\wedge G)$
is achieved either for $F=G=\Q$ and equals $d_A(\Q)+d_B(\Q)=d_C(\Q)$,
or is achieved
for $F=\Z/\p^\infty$, $G=\Z/\p^\infty$ and is equal to $\alpha(F)+\beta(G)+1$,
or is achieved either for $F=\Z/\p$, $G=\Z/\p$ and is equal to
$\alpha(F)+\beta(G)=d_A(\Z/\p)+d_B(\Z/\p)$.
From what we know about $d_C(\Z/{\p^{\infty}})$ we conclude
$d_C(\Z_{(\p)})
=\min(d_C(\Q),d_C(\Z/{\p^{\infty}})-1).$
 \end{pf}


\begin{Def} \label{XX4.40} Given two Bockstein functions
$\alpha:\B\to\Z\cup\{\pm\infty\}$
and $\beta:\B\to\Z$ define their {\bf smash product} $\alpha\wedge \beta$
as the extension
 function
of $GG(\alpha)\wedge GG(\beta)$.

\end{Def}

The next result is a corollary to \ref{XX4.39} and it shows that our definition
of the smash product of Bockstein functions coincides with the one given in
\cite{D-D$_1$}.


\begin{Thm} \label{XX4.41} Given two Bockstein functions
$\alpha:\B\to\Z\cup\{\pm\infty\}$
and $\beta:\B\to\Z\cup\{\pm\infty\}$ their smash product $\alpha\wedge
\beta$ is given by
the following formulae:
\par{1.} Let $\gamma=\alpha+\beta$,
\par{2.} If $H\ne \Z/{\p^{\infty}}$ and $H\ne \Z_{(\p)}$
for all $\p$, then
$$(\alpha\wedge \beta)(H)=\gamma(H),$$
\par{3.} If one of $\alpha,\beta$ is $\p$-regular for some $\p$, then
$$(\alpha\wedge \beta)(H)=\gamma(H)$$
for $H=\Z/{\p^{\infty}},\Z_{(\p)}$.

\par{4.} If both $\alpha,\beta$ are $\p$-singular for some $\p$, then
$$(\alpha\wedge
\beta)(\Z/{\p^{\infty}})=\min(\gamma(\Z/{\p}),\gamma(\Z/{\p^{\infty}})+1),$$
$$(\alpha\wedge \beta)(\Z_{(\p)})
=\max(\gamma(\Q),(\alpha\wedge \beta)(\Z/{\p^{\infty}})+1).$$

\end{Thm}


\begin{Def} \label{XX4.42} Given a Bockstein function
$\beta:\B\to\Z\cup\{\pm\infty\}$
 define its {\bf dual} $\beta^*$ as
the extension function of $(GG(\beta))^*$.

\end{Def}


\begin{Thm} \label{XX4.43} Given a Bockstein function
$\beta:\B\to\Z\cup\{\pm\infty\}$
its dual $\beta^*$ is given by the following formulae:
\par{1.} $\beta^*(H)=-\beta(H)$ for $H=\Q,\Z/\p$.
\par{2.} If $\beta$ is $\p$-regular for some $\p$, then
$\beta^*(H)=-\beta(H)$ for $H=\Z/\p^\infty,\Z_{(p)}$.
\par{3.} If $\beta$ is $\p$-singular for some $\p$, then
$\beta^*(\Z/\p^\infty)=-\beta(\Z/\p^\infty)-1$
and
$\beta^*(\Z_{(p)})=\max(-\beta(\Q),-\beta(\Z/\p^\infty))$.

\end{Thm}

\begin{pf}  Let $A=GG(\beta)$ and $B=A^*$. By \ref{XX4.25},
$d_B=-e_A=-\beta$. Using \ref{XX4.38} one arrives at  1)-3).
 \end{pf}


\begin{Thm} \label{XX4.44} Let $\cal{BF}$ be the family of all Bockstein
functions.
If $\beta$ is a Bockstein function, then
 $$\beta^*=\inf\{\alpha\in\cal{BF}\mid \alpha\wedge\beta\ge 0\}.$$

\end{Thm}

\begin{pf}  Notice that $\beta^*\wedge\beta\ge 0$.
Suppose $\alpha\wedge\beta\ge 0$ and let $A=GG(\alpha)$,
$B=GG(\beta^*)$. Since $\dim_\Z(A\wedge B^*)\ge 0$,
$\dim(B)\leq \dim(A)$ (see \ref{XX4.21})
and $\beta^*=e_B\leq e_A=\alpha$ (see \ref{XX4.32} and \ref{XX4.34}).
 \end{pf}


\begin{Def} \label{XX4.45}
Given two Bockstein
functions $\alpha$ and $\beta$, one defines their {\bf sum--product}
 $\alpha\lbrack +\rbrack \beta$ as follows:
\par{1.} $(\alpha\lbrack +\rbrack \beta)(\Z/\p)=\alpha(\Z/\p)+\beta(\Z/\p)$,
\par{2.} $(\alpha\lbrack +\rbrack \beta)(\Q)=\alpha(\Q)+\beta(\Q)$,
\par{3.} $(\alpha\lbrack +\rbrack \beta)(\Z/{\p^{\infty}})=
\max\{\alpha(\Z/{\p^{\infty}})+\beta(\Z/{\p^{\infty}}),
\alpha(\Z/\p)+\beta(\Z/\p)-1\}$,
\par{4.} $(\alpha\lbrack +\rbrack
\beta)(\Z_{(\p)})=\alpha(\Z_{(\p)})+\beta(\Z_{(\p)})$ if
$\alpha(\Z_{(\p)})=\alpha(\Z/{\p^{\infty}})$ or
$\beta(\Z_{(\p)})=\beta(\Z/{\p^{\infty}})$,
\par $(\alpha\lbrack +\rbrack \beta)(\Z_{(\p)})=
\max\{(\alpha\lbrack +\rbrack
\beta)(\Z/{\p^{\infty}})+1,\alpha(\Q)+\beta(\Q)\}$
if $\alpha(\Z_{(\p)})>\alpha(\Z/{\p^{\infty}})$ and
$\beta(\Z_{(\p)})>\beta(\Z/{\p^{\infty}})$.
 
\end{Def}

\begin{Rem}  $\alpha\lbrack +\rbrack \beta$ was introduced in
\cite{K} under the name of the product of two Bockstein functions
and is denoted there by $\alpha\times \beta$. A.Dranishnikov
\cite{D$_7$} realized that $\alpha\lbrack +\rbrack \beta$
is much closer related to the sum $\alpha + \beta$ and that
there is another operation, denoted by $\alpha\lbrack \times\rbrack \beta$,
which should play the role of a product.
 
\end{Rem}


\begin{Thm} \label{XX4.46} $\alpha\lbrack +\rbrack
\beta=(\alpha^*\wedge\beta^*)^*$
for every two Bockstein functions $\alpha$ and $\beta$.

\end{Thm}

\begin{pf}  Let $\gamma=\alpha\lbrack +\rbrack \beta$
and $\nu=(\alpha^*\wedge\beta^*)^*$.
Applying \ref{XX4.41} and \ref{XX4.43} we get $\gamma(H)=\nu(H)$ if $H=\Z/\p,\Q$.
\par
If both $\alpha$ and $\beta$ are $\p$-regular, then
$\gamma$ and $\nu$ are $\p$-regular and
$\gamma(H)=\nu(H)=\alpha(H)+\beta(H)$ for $H=\Z_{(p)},\Z/\p^\infty$.
\par
If both $\alpha$ and $\beta$ are $\p$-singular, then both $\alpha^*$ and
$\beta^*$ are $\p$-singular.
By \ref{XX4.43}, $\alpha^*(\Z/\p^\infty)=-\alpha(\Z/\p^\infty)-1$
and $\beta^*(\Z/\p^\infty)=-\beta(\Z/\p^\infty)-1$.
By \ref{XX4.41}, $(\alpha^*\wedge\beta^*)(\Z/\p^\infty)=
\min(-\alpha(\Z/\p)-\beta(\Z/\p),-\alpha(\Z/\p^\infty)-\beta(\Z/\p^\infty)-1)$.
Applying \ref{XX4.41} again, we get $\nu(\Z/\p^\infty)=
\max(\alpha(\Z/\p)+\beta(\Z/\p),\alpha(\Z/\p^\infty)+\beta(\Z/\p^\infty)+1)-1=
\gamma(\Z/\p^\infty)$. Now, both $\gamma$ and $\nu$ are $\p$-singular
and agree on $\Q$ and $\Z/\p^\infty$ which implies that
they agree on $\Z_{(p)}$.
\par
Assume $\alpha$ is $\p$-regular and $\beta$ is $\p$-singular.
Now, $\alpha^*$ is $\p$-regular and $\beta^*$ is $\p$-singular.
By \ref{XX4.43}, $\beta^*(\Z/\p^\infty)=-\beta(\Z/\p^\infty)-1$,
$\beta^*(\Z_{(p)})=\max(-\beta(\Q),-\beta(\Z/\p^\infty))$, and
$\alpha^*(\Z/\p^\infty)=-\alpha(\Z/\p^\infty)=\alpha^*(\Z_{(p)})$.
Applying \ref{XX4.41} we get
$(\alpha^*\wedge\beta^*)(H)=\alpha^*(H)+\beta^*(H)$
for $H=\Z_{(p)},\Z/\p^\infty$. By \ref{XX4.43},
$\nu(\Z/\p^\infty)=-\alpha^*(\Z/\p^\infty)-\beta^*(\Z/\p^\infty)-1
=\alpha(\Z/\p^\infty)+\beta(\Z/\p^\infty)=\gamma(\Z/\p^\infty)$.
Again, both $\gamma$ and $\nu$ are $\p$-singular
and agree on $\Q$ and $\Z/\p^\infty$ which implies that
they agree on $\Z_{(p)}$.
 \end{pf}


\begin{Thm} \label{XX4.47} For any dimension-like function
 $d:\cal{DGG}\to\Z\cup\{\pm\infty\}$ there is a graded group $A$
such that $d(X)=\dim_\Z(A\wedge X)$.

\end{Thm}

\begin{pf}  $e=-d$ is an extension function and there is a graded group $B$
so that $e(H)=e_B(H)$ for each $H\in\B$. Let $A=B^*$. Now,
$d_A(H)=-e_B(H)=d(H)$ for each $H\in\B$ which means $d_A=d$.
 \end{pf}


\begin{Thm} \label{XX4.48} A function $\alpha:\B\to\Z\cup\{\pm\infty\}$
extends to an extension function $e:\cal{DGG}\to\Z\cup\{\pm\infty\}$
if and only if $\alpha$ is a Bockstein function.

\end{Thm}

\begin{pf}  By \ref{XX4.31}, a restriction
$e\mid \B$ of an extension function $e:\cal{DGG}\to\Z\cup\{\pm\infty\}$
is a Bockstein function.
Suppose $\beta:\B\to \Z\cup\{\pm\infty\}$ is a Bockstein function.
Let $A=GG(\beta)$ be the associated graded group. \ref{XX4.34} says that
$\beta=e_a\mid \B$ and $e_A$ is an extension function by \ref{XX4.27}.
 \end{pf}


\begin{Thm} \label{XX4.49} Any set $\{A_t\}_{t\in T}$ of $\cal{DGG}$
has infimum $\bigoplus\limits_{t\in T}A_t$
and supremum $(\bigoplus\limits_{t\in T}A_t^*)^*$.

\end{Thm}

\begin{pf}  Clearly, $\dim(\bigoplus\limits_{t\in T}A_t)\leq\dim(A_t)$
for each $t\in T$. If $\dim(B)\leq \dim(A_t)$ for each $t\in T$, then
it is easy to show that $\dim(B)\leq\dim(\bigoplus\limits_{t\in T}A_t)$.
\par Dually, $\dim(\bigoplus\limits_{t\in T}A_t^*)\leq\dim(A_t^*)$
for each $t\in T$ which implies
$\dim((\bigoplus\limits_{t\in T}A_t^*)^*)\ge\dim(A_t)$
for each $t\in T$.
If $\dim(B)\ge \dim(A_t)$ for each $t\in T$, then
$\dim(B^*)\leq \dim(A_t^*)$ and
 $\dim(B^*)\leq\dim(\bigoplus\limits_{t\in T}A_t^*)$.
Therefore, $\dim(B)\ge\dim((\bigoplus\limits_{t\in T}A_t^*)^*)$.
 \end{pf}


\begin{Cor} \label{XX4.50} $\cal{DGG}$ is a lattice. Moreover,
$\min(A,B)=A\oplus B$ and $\max(A,B)=(A^*\oplus B^*)^*$.

\end{Cor}

It would be interesting to characterize self-dual graded groups.


\begin{Prop} \label{XX4.51} Suppose $A$ is a graded group. $\dim(A)=\dim(A^*)$
if and only if $\dim(A)=\dim(\Z)$.

\end{Prop}

\begin{pf}  Notice that $e_\Z=0$ which implies $e_{\Z^*}=0$ by \ref{XX4.43}.
Thus, $\dim(\Z)=\dim(\Z^*)$.
If $\dim(A)=\dim(A^*)$, then setting $\beta=e_A$ one gets
$\beta=\beta^*$ which is possible only if $\beta=0$ (see \ref{XX4.43}),
i.e. $\dim(A)=\dim(\Z)$ (see \ref{XX4.32}).
 \end{pf}

\section{Applications}

\begin{Def} \label{XXa5.3} Suppose $K$ is a pointed CW complex. Given an Abelian group $G$
define define $e^K(G)$ as $e_A(G)$, where $A=H_*(K)$.
\end{Def} 

Functions $e^K$ were first introduced in \cite{D-D$_1$}
to study extension properties of CW complexes via so-called Dual Bockstein Algebra.
The idea was to dualize the approach of \cite{K}
where Bockstein functions $d_X$ were used to investigate cohomological
dimension of compact spaces $X$. Thus, one has
two families of Bockstein functions $d_X,e^K:\B\to \Z\cup\{\pm\infty\}$.

\begin{Cor} \label{XXa5.5} Suppose $X$ is a pointed compact
space and $K$ is a pointed CW complex.
The following conditions are equivalent:
\par{1.} $SP(K)\in AE(X)$.
\par{2.} $d_X\leq e^K$.

\end{Cor}

\begin{pf}  By  \ref{XX5.5} $SP(K)\in AE(X)$ if and only if
$\cal H^{-\ast}(X)\wedge H_*(K)$ is non-negative. Therefore
1) is equivalent to $(\cal H^{-\ast}(X))^*\leq H_*(K)$
which is equivalent to 2).
 \end{pf}

The following result was proved in \cite{D-D$_1$} (see Theorem 5.20)
for countable CW complexes.


\begin{Cor} \label{XX5.6} If $K$, $L$ are pointed CW complexes,
then $e^{K\wedge L}=e^K\wedge e^L$.

\end{Cor}

\begin{pf}  $e^K\wedge e^L$ was defined in \ref{XX4.40}
as the extension function of $GG(e^K)\wedge GG(e^L)$.
By \ref{XX4.33} and \ref{XX4.35}, $GG(e^M)$ has the same dimension as $H_*(M)$
for each CW complex $M$. Thus, $e^K\wedge e^L$ is
 the extension function of $H_*(K)\wedge H_*(L)$
which is $e^{K\wedge L}$ in view of \ref{XXeh.1}.
 \end{pf}


\begin{Def} \label{XX5.8} Given a function $\alpha:\B\to [1,\infty]$
define the Moore space $M(\alpha)$ as the wedge of those $M(H,\alpha(H))$
for which $\alpha(H)$ is finite.

\end{Def}

The following result was proved in \cite{Dy$_6$} using different methods.

\begin{Thm} \label{XX5.9} Suppose $X$ is a pointed compact space. There is a minimum
$\{SP(K)\mid SP(K)\in AE(X)\}$ called the
{\bf cohomological dimension} of $X$. That minimum is represented by a
countable CW complex.

\end{Thm}

\begin{pf}  If $X$ is totally disconnected,
then $S^0\in AE(X)$ and $SP(S^0)\sim S^0$ is the minimum
of all CW complexes $K$ so that $K\in AE(X)$.
Suppose $X$ is not totally disconnected.
Now, $A=(\cal H^{-\ast}(X))^*$ is positive and we may consider $K=M(A')$,
where $A'$ is built of Bockstein groups and $\dim(A')=\dim(A)$.
Notice that $d_X=e^K$ which implies $SP(K)\in AE(X)$ in view of \ref{XXa5.5}.
If $SP(L)\in AE(X)$ for some $L$,
then $e_K=d_X\leq e_L$ and $H_\ast(K)\leq H_\ast(L)$ which implies
$SP(K)\leq SP(L)$ (see \ref{XXa5.5}).
 \end{pf}


\begin{Thm} [First Bockstein Theorem]\label{XX5.11}  For every compact space $X$
and every Abelian group $G\ne 0$
\par
\centerline{$\dim_G(X)=\sup\{\dim_H(X)\mid H\in\sigma(G)\}.$}

\end{Thm}

\begin{pf}  Let $F=\bigoplus\limits_{H\in\sigma(G)}H$.
By \ref{XX4.16}, $\dim_G(X)=\dim_F(X)$ and it is clear
that $\dim_F(X)=\sup\{\dim_H(X)\mid H\in\sigma(G)\}.$
 \end{pf}


\begin{Thm} [Second Bockstein Theorem]\label{XX5.12} 
Suppose $X$ and $Y$
are pointed compact spaces. Then,
\par{1.} $\dim_{\Z/\p}(X\wedge Y)=\dim_{\Z/\p}X+\dim_{\Z/\p}Y$,
\par{2.} $\dim_{\Q}(X\wedge Y)=\dim_{\Q}X+\dim_{\Q}Y$,
\par{3.} $\dim_{\Z/{\p^{\infty}}}(X\wedge Y)=
\max\{\dim_{\Z/{\p^{\infty}}}X+\dim_{\Z/{\p^{\infty}}}Y,
\dim_{\Z/\p}X+\dim_{\Z/\p}Y-1\}$,
\par{4.} $\dim_{\Z_{(\p)}}(X\wedge Y)=\dim_{\Z_{(\p)}}X+\dim_{\Z_{(\p)}}Y$
if
$\dim_{\Z_{(\p)}}X=\dim_{\Z/{\p^{\infty}}}X$ or
$\dim_{\Z_{(\p)}}Y=\dim_{\Z/{\p^{\infty}}}Y$,
\par $\dim_{\Z_{(\p)}}(X\wedge Y)=
\max\{\dim_{\Z/{\p^{\infty}}}(X\times Y)+1,\dim_{\Q}X+\dim_{\Q}Y\}$
if $\dim_{\Z_{(\p)}}X>\dim_{\Z/{\p^{\infty}}}X$ and
$\dim_{\Z_{(\p)}}Y>\dim_{\Z/{\p^{\infty}}}Y.$

\end{Thm}

\begin{pf} Let $A=\cal H^{-\ast}(X)$,
$B=\cal H^{-\ast}(Y)$, and $C=\cal H^{-\ast}(X\wedge Y)$.
Notice that $\dim(C)=\dim(A\wedge B)=\dim((A^*)^*\wedge (B^*)^*)$ (see \ref{XX5.2})
which implies $e_{C^*}=e_{A^*}\lbrack +\rbrack e_{B^*}$.
Since $d_X=-d_A=e_{A^*}$,
$d_Y=-d_B=e_{B^*}$,
and $d_{X\wedge Y}=-d_C=e_{C^*}$, \ref{XX5.12} follows from \ref{XX4.46}.
 \end{pf}


\begin{Thm} [Dranishnikov Realization Theorem]\label{XX5.13}  Suppose $n\ge 1$
and $\alpha:\B\to [1,n]$ is a Bockstein function.
There is a compactum $X\subset I^{n+2}$ such that
$d_X=\alpha$.

\end{Thm}

\begin{pf}  Let $A=GG(\alpha)$, $B=\Sigma^{n+1}(A^*)$,
and $\beta=e_B$. Notice that $\beta\ge 1$. Indeed,
using \ref{XX4.43} one gets $(\alpha)^*\ge -n$ and it is clear
that $\beta=(\alpha)^*+n+1$ which implies $\alpha\wedge\beta\ge n+1$.
Let $K=M(\alpha)$ and $L=M(\beta)$. Notice that
the extension function of $K\wedge L$ is $\alpha\wedge\beta\ge n+1$
which means that the extension function of $K*L$ is at least $n+2$,
i.e. $K*L\in AE(I^{n+2})$. Split $I^{n+2}$ as $X'\cup Y$
so that $X'$ is $\sigma$-compact, $K\in AE(X')$, and $L\in AE(Y)$.
Replace $X'$ by the compact space $X$ as follows: if $X'$ is the union
of its compact subspaces $X_n$, $n\ge 1$, then $X$ is the compact wedge
of all $X_n$. In particular, $X'$ and $X$ have the same
extension dimension (that means $K\in AE(X)$ is equivalent to $K\in AE(X')$
for all CW complexes $K$) and
 $d_X\leq \alpha$. Let $d_X=\gamma$ and choose $M$ with extension function
$\gamma$. Now, $M*L\in AE(I^{n+2})$, i.e. $\gamma\wedge\beta\ge n+1$
and $\gamma\ge (\beta)^*+n+1=\alpha$. This proves
$d_X= \alpha$.
 \end{pf}

\begin{Cor} \label{XX5.13.5} For every Bockstein function
$\alpha:\B\to [1,\infty]$
there is a pointed compactum $X$ such that
$d_X=\alpha$.
\end{Cor}

\begin{pf} For each $n\ge 1$ put $\alpha _n=\min(\alpha ,n)$
and pick a pointed compactum $X(n)$ so that $d_{X(n)}=\alpha _n$.
The compact wedge $X$ of all $X_n$ is the required pointed compactum.
 \end{pf}


\begin{Thm} [Test Spaces Theorem]\label{XX5.14}  For each $n>1$ and each Abelian
group $G$ there is a pointed compactum $T\subset I^{n+2}$ such that
$\dim_\Z(X\wedge T)=\dim_G(X)+n$ for every pointed compact space $X$
satisfying $\dim_\Z(X)<\dim_G(X)+n$.

\end{Thm}

\begin{pf}  Consider the graded group $B=\Sigma^{-n}(G)\oplus\Sigma^{-1}(\Z)$.
We need to find a pointed compactum $T$ so that
$\dim(B)=\dim(\cal H^{-\ast}(T))$. Consider $A=B^*$.
Notice that $-n\leq d_B\leq -1$
which implies $1\leq e_A\leq n$ (see \ref{XX4.25}).
In particular, $A$ is positive.
Choose $T\subset I^{n+2}$
with $d_T=e_A$ (see \ref{XX5.13}) which implies
$\dim(\cal H^{-\ast}(T))=\dim(B)$.
If $\dim_G(X)=k$
and $\dim_\Z(X)<k+n$,
then $C=\cal H^{-\ast}(X;G)$ satisfies $C(-k)\ne 0$
and $C(i)=0$ for $i<-k$.
Look at $D=\cal H^{-\ast}(X)\wedge \cal H^{-\ast}(T)$ and notice
that $D(-k-n)\ne 0$ and $D(i)=0$ for $i<-k-n$.
That means $\dim_\Z(X\wedge T)=k+n$ by \ref{XX5.2}.
 \end{pf}

\begin{Rem}  Notice that \ref{XX5.14} is slightly more general than
6.1 of \cite{D$_8$} in the sense that it deals with
integral dimension instead of the covering dimension.
Also notice that 6.1 of \cite{D$_8$} is a significant
improvement of the original Test Spaces Theorem
(see Theorem 12 in \cite{K}).
 
\end{Rem}

\section{Embeddings of algebras in $\cal{DGG}$}

In this section we will discuss embeddings of geometrically defined
algebras in the Dimension Algebra of Graded Groups.

\begin{Prop}\label{XXr.0} Let $\cal{CW}$ be the subalgebra
of the Standard Algebra consisting of all pointed CW complexes.
Assigning $\dim(H_\ast(K))$ to $K$
is a homomorphism from $\cal{CW}$ to $\cal{DGG}$
whose image consists of $\dim(A)$ such that $A$ is a non-negative
graded group and $A(0)$ is free Abelian.
\end{Prop}
\begin{pf} Clearly, $H_\ast(K\vee L)\equiv H_\ast(K)\oplus H_\ast(L)$.
Also $H_\ast(K\wedge L)\equiv H_\ast(K)\wedge H_\ast (L)$
by \ref{XXeh.1}.
\par Obviously, $\dim(H_\ast(K))$ is a non-negative graded group
whose 0-th term is free Abelian. Suppose $A$ is a non-negative
graded group and $A(0)$ is free Abelian.
\par Case 1. $A(0)=0$. In this case one can consider the Moore
space $M_n=M(A_n,n)$ (a pointed CW complex with only one non-zero homology
and whose $n$-th homology is $A_n$). The wedge $M$ of all $M_n$, $n\ge 1$,
is a CW complex such that $H_\ast(M)=A$.
\par Case 2. $A(0)\ne 0$. In this case it is obvious that $\dim(A)=\dim(Z)$
and $M=S^0$ satisfies $\dim(H_\ast(M))=\dim(A)$.
\end{pf} 

\begin{Prop}\label{XXr.00} Let $\cal{COMP}$ be the subalgebra
of the Standard Algebra consisting of all pointed compact spaces.
Assigning $\dim(\cal H^{-\ast}(X))$ to $X$
is a homomorphism from $\cal{COMP}$ to $\cal{DGG}$
whose image consists of $\dim(0)$ and $\dim(A^\ast)$
such that $A$ is a non-negative
graded group and $A(0)$ is free Abelian.
\end{Prop}
\begin{pf} Clearly, $\cal H^{-\ast}(X\vee Y)\equiv \cal H^{-\ast}(X)\oplus \cal H^{-\ast}(Y)$.
Also $\cal H^{-\ast}(X\wedge Y)\equiv \cal H^{-\ast}(X)\wedge \cal H^{-\ast} (Y)$
by \ref{XX5.2}.
\par If $X$ is a pointed point, then $\cal H^{-\ast}(X)=0$.
If $X$ contains at least two points, then put  $B=(H^{-\ast}(X))^\ast$.
Notice that $e_B=d_X$, so it is either positive or identically $0$
(see \ref{XXeh.2}). Since $C=GG(e_B)$ satisfies $\dim(C)=\dim(B)$
by \ref{XX4.35}, $C$ is positive if $e_B > 0$ or $\dim(C)=\dim(\Z)$ if
$e_B\equiv 0$. Putting $A=C$ if if $e_B > 0$ and $A=\Z$ if
$e_B\equiv 0$ one gets a non-negative graded group $A$
such that $A(0)$ is free and $\dim(A^\ast)=\dim(\cal H^{-\ast}(X))$.
\par Conversely, if $A$ is a non-negative graded group such that $A(0)$ 
is Abelian, then $e_A$ is a Bockstein function
which is either positive or identically $0$ if $A(0)\ne 0$.
Therefore, by  \ref{XX5.13.5}, there is a pointed compactum $X$
such that $d_X=e_A$. That amounts to $\dim(\cal H^{-\ast}(X))=\dim(A^\ast)$
by \ref{XX4.25}.
\end{pf}

Let us give sufficient and necessary conditions
for two pointed CW complexes in $\cal{CW}$ to be mapped to the same element of
$\cal{DGG}$. The result below was proved in \cite{Dy$_6$}
for countable CW complexes using different methods.

\begin{Prop}\label{XXr.0.5} If $K$ and $L$ are pointed CW complexes, then
the following conditions are equivalent:
\par 1. $\dim(H_\ast(K))=\dim(H_\ast(L))$.
\par 2. $cin(K\wedge M)=cin(L\wedge M)$ for all pointed CW complexes $M$.
\par 3. There is $n\ge 2$ such that $\Sigma^n(K)\sim_X \Sigma^n(L)$ for all finite-dimensional
compacta $X$.
\par 4. $SP(K)\sim_X SP(L)$ for all compact spaces $X$.
\end{Prop}
\begin{pf} 1)$\implies$2). Notice that $cin(P\wedge M)=\dim_\Z(H_\ast(P)\wedge H_\ast(M))$
for all pointed CW complexes $M$ and $P$. Using  \ref{XXeh.1}, 1)$\implies$2) follows.
\par 2)$\implies$1). Given an abelian group $G$ consider $M=M(G,1)$, the Moore space.
Thus, $H_\ast(P)=\Sigma(G)$ and $cin(P\wedge M)=\dim_\Z(H_\ast(P)\wedge \Sigma(G))=
1+\dim_\Z(H_\ast(P)\wedge G)=1+\dim_G(H_\ast(P))$.
Thus $1+\dim_G(H_\ast(K))=cin(K\wedge M)=cin(L\wedge M)=1+\dim_G(H_\ast(L))$
which means $\dim_G(H_\ast(K))=\dim_G(H_\ast(L))$ for all $G$,
i.e. $\dim(H_\ast(K))=\dim(H_\ast(L))$.
\par 1)$\iff$3) and 1)$\iff$4) follow from \ref{XXeh.9} and \ref{XXeh.7}.
\end{pf}

The following two results are immediate corollaries of
\ref{XXr.0.5} and \ref{XXr.0}.

\begin{Thm}\label{XXr.1} Assigning $\dim(H_\ast(K))$ to
the equivalence class of a pointed CW complex $K$
induces an embedding of Shchepin Algebra into $\cal{DGG}$.
\end{Thm}

\begin{Thm}\label{XXr.2} Assigning $\dim(H_\ast(K))$ to
the equivalence class of a pointed CW complex $K$
induces an embedding of Dranishnikov-Dydak Algebra into $\cal{DGG}$.
\end{Thm}

\begin{Cor}\label{XXr.3} Dranishnikov-Dydak Algebra and
Shchepin Algebra are identical.
\end{Cor}
\begin{pf} In view of the above results one has
the natural embedding of Dranishnikov-Dydak Algebra
to Shchepin Algebra. Notice that \ref{XX4.19} implies that
for any pointed CW complex $K$ there is a pointed countable
CW complex $K'$ with $\dim(H_\ast(K))=\dim(H_\ast(K'))$.
Therefore, the embedding of algebras is surjective.
\end{pf}

\begin{Thm}\label{XXr.4} Assigning $\dim(\cal H^{-\ast}(X))$
to $X$
induces an embedding of Kuzminov Algebra into $\cal{DGG}$.
\end{Thm}
\begin{pf} Suppose $X$ and $Y$ are two pointed compact spaces
of finite dimension such that
$\dim(X\wedge T)=\dim(Y\wedge T)$ for all pointed compact spaces $T$
of finite dimension.
To show $\dim(\cal H^{-\ast}(X))=\dim(\cal H^{-\ast}(Y))$
it suffices to prove $\dim_G(X)=\dim_G(Y)$ for all Abelian groups $G$.
Pick $n > \dim(X),\dim(Y)$ and use  \ref{XX5.14}
to produce $T$ such that $\dim(X\wedge T)=\dim_G(X)+n$
and $\dim(Y\wedge T)=\dim_G(Y)+n$.
Hence $\dim_G(X)=\dim_G(Y)$.
\par Conversely, if $\dim(\cal H^{-\ast}(X))=\dim(\cal H^{-\ast}(Y))$,
i.e.  $\dim_G(X)=\dim_G(Y)$ for all Abelian groups $G$,
then \ref{XX5.12} implies $\dim_H(X\wedge T)=\dim_H(Y\wedge T)$
for all Bockstein groups $H$.
Applying \ref{XX5.11} one gets $\dim_\Z(X\wedge T)=\dim_\Z(Y\wedge T)$
which implies $\dim(X\wedge T)=\dim(Y\wedge T)$ if $T$ is finite-dimensional.
\end{pf} 

Proposition 3.7 of \cite{D-D$_1$} gives the following improvement
of \ref{XXea.2}. In contrast to \ref{XXea.2} that result
requires algebra which is quite easy in our setting.

\begin{Lem}\label{XXr.5}
Suppose $K_1,K_2,L_1,L_2$ are pointed countable CW complexes.
If $K_i\sim_X L_i$ for $i=1,2$ and for all finite-dimensional compacta $X$,
then $K_1\wedge L_1\sim_X K_2\wedge L_2$ for all finite-dimensional compacta $X$.
\end{Lem}
\begin{pf} As in Proposition 3.7 of \cite{D-D$_1$}
the only non-geometric case is that of connected CW complexes.
In that case $K_1\wedge L_1\sim_X K_2\wedge L_2$ for all finite-dimensional compacta $X$
if $H_\ast(K_1\wedge L_1)$ and $H_\ast(K_2\wedge L_2)$
represent the same element of $\cal{DGG}$.
However, that is true as \ref{XXea.2}
and  \ref{XXeh.9} imply that their suspensions represent the same element of $\cal{DGG}$
(see \ref{XXeh.9}).
\end{pf}

Therefore the following definition makes sense.
\begin{Def}\label{XXr.6}
The {\bf Unstable Extension Algebra} is the quotient
of the Standard Algebra of pointed countable CW complexes
under the equivalence relation $K\sim_X L$ for all
finite-dimensional compacta.
\end{Def}

Notice that the above algebra is not identical with the Dranishnikov-Dydak
Algebra. Indeed, any $M=M(G,1)$, where $G\ne 0$ is perfect,
is not equivalent to $I$ as $M\notin AE(I^2)$ and $I\in AE(I^2)$.
On the other hand $\Sigma(M)$ and $\Sigma(I)$ are both contractible.

\medskip
\medskip
\medskip
\medskip
Jerzy Dydak,
Math Dept, University of Tennessee, Knoxville, TN 37996-1300, USA,
E-mail address: dydak@@math.utk.edu

\end{document}